# Mathematical model of COVID-19 with imperfect vaccine and virus mutation

Ceren Gürbüz Can[1], Sebaheddin Şevgin[1]

*Correspondence:
ceren.gurbuz@hotmail.com

[1] Department of Mathematics, Faculty of Sciences, Van Yüzüncü Yıl University. Van, Turkey.

**Abstract**

This study examines the effect of a partially protective vaccine on COVID-19 infection with the original and mutant virus with the help of a deterministic mathematical model developed. The model we developed consists of five compartments and thirteen parameters. The model consists of $S$ (susceptible), $V$ (vaccinated), $I_1$ (infected with original virus), $I_2$ (infected with mutant virus) and $R$ (recover) subcompartments. With the established model, imperfect vaccines and mutant virus for the COVID-19 pandemic in Turkey were examined. We examined the effect of both artificial active immunity (vaccinated) and natural active immunity (passing disease) in the model. Since it is known that the recovery and death rates of the original virus and the mutant virus are different in COVID-19, we considered it in the study. We performed local stability and global stability analysis by calculating the disease-free equilibrium point and endemic equilibrium point of the model. We also obtained the basic reproduction number with the help of the next generation matrix method. We estimated three model parameters with parameter estimation and identified model-sensitive parameters with local sensitivity analysis. Next, we showed the existence of a backward bifurcation using the Castillo-Chavez and Song Bifurcation Theorem. Additionally, to illustrate the underlying fundamental mechanisms of the model dynamics and to support the analytical findings, we presented three distinct simulations.

## 1 Introduction

Infectious diseases have found opportunities to spread more rapidly and easily with the changing social living conditions. Cholera, Dengue fever, Ebola, and the COVID-19 pandemic, which we are currently fighting, are among the most significant outbreaks identified to date. In addition to their deadly consequences, infectious diseases have had severe demographic, economic, and social impacts on a global scale. Particularly, COVID-19 has become an important example in demonstrating how profoundly infectious diseases can affect human life, as it has impacted the entire world.

The COVID-19 pandemic is a viral outbreak that originated in Wuhan, the capital of Hubei Province, China, on November 17, 2019. The virus, which is transmissible from person to person, saw a significant increase in its transmission rate by mid-January 2020. As the pandemic progressed, cases were reported in various countries across Europe, North America, and the Asia-Pacific region. The first confirmed case of COVID-19 in Turkey was announced by the Ministry of Health on March 11, 2020. The first death related to the virus in the country occurred on March 15, 2020. The vaccination campaign began on January 14, 2021. The pandemic led to the implementation of radical decisions, which had significant social, economic, political, administrative, legal, military, religious, and cultural effects and consequences in Turkey. Hundreds of settlements, villages, and towns were placed under quarantine as part of COVID-19 containment measures [6].

Mathematical epidemiology, which scientifically examines the process of mathematical modeling of infectious diseases, focuses on viral and bacterial diseases. These diseases (viral or bacterial), which are the subject of mathematical epidemiology, have been modeled using differential equation systems with the help of various subdivisions since 1927. In these models, the total population is categorized based on exposure to the pathogen causing the disease, treatment, or immunity status. Thus, the population is divided into sub-compartments such as susceptible, vaccinated, infected, and recovered.

Theoretical papers on disease models that fall within the scope of epidemiology have become important tools in analyzing the spread and control of diseases. Additionally, mathematical models are intended to be used for comparing, planning, implementing, evaluating, and optimizing various detection, prevention, treatment, and control programs.

In the literature, mathematical models used to study the dynamics of infectious diseases are classified as dynamic and statistical models. Dynamic models investigate the evolution of disease processes over time. These models include deterministic, stochastic, network, and agent-based models. Deterministic models, one of which is examined in our study, are also known as compartmental models as they are formulated by dividing the total population into subpopulations. These models do not account for the effects of randomness in disease dynamics. They are suitable for modeling large populations and use differential equations for continuous-time models and difference equations for discrete-time models [17].

One of the earliest mathematical models in epidemiology was related to the effects of vaccination. In 1766, Swiss mathematician Daniel Bernoulli published a study on the estimated impact of smallpox vaccination on the life expectancy of the vaccinated population [2].

A significant contribution to mathematical epidemic theory was made through the work of W.O. Kermack and A.G. McKendrick. The basic compartmental models that describe the transmission mechanisms of infectious diseases are presented in a series of three papers published by W.O. Kermack and A.G. McKendrick in 1927, 1932, and 1933 [18].

In [15], T. D. Frank analyzed the dynamics of COVID-19 infections from the perspective of bifurcation theory by modeling the $R$ (recovered) class in the basic SEIR compartmental model as five separate subpopulations. In [25], the $D$ (death) subpopulation was added to the basic $SEIR$ model, and analytical solutions were examined. In [12], possible equilibrium points were calculated using the $SIQ$ compartmental model, and stability analysis was conducted using Lyapunov theory. A. E. Matouk, in [22], examined the effects of two therapeutic drugs on the $SIMDR$ model, identifying eight equilibrium points and performing a stability analysis. A. B. Gumel et al., in [17], modeled the $I$ (infectious) class of the basic $SEIR$ model as three distinct subpopulations, analyzing the equilibrium points and stability. In [8], asymptomatic infectious and hospitalized infectious subpopulations were added to the standard infectious subpopulation of the classical $SEIR$ model, and local asymptotic stability properties were explored based on the characteristics of the next-generation matrix, the value of the basic reproduction number, and equilibrium states.

In [5], the pandemic dynamics were analyzed by comparing the data from Italy, which experienced the COVID-19 pandemic early, and Turkey, which faced it later, using a five-compartment model. In [1], the effects of social distancing on the dynamics of the pandemic were observed and a model developed using information from the Hubei outbreak was applied to analyze the dynamics of the pandemic in Turkey. In [24], a mathematical prediction model was used to estimate confirmed cases and deaths caused by COVID-19 in Iran and Turkey. In [16], the COVID-19 pandemic in Turkey was modeled using a five-compartment approach, and both local and global stability analyses were conducted.

In this study, a deterministic compartmental model was constructed and analyzed to evaluate the effects of the original and mutant COVID-19 viruses and partial protective vaccination on pandemic eradication. The developed model included vaccinated individuals as an additional compartment. However, issues related to vaccination, such as some vaccines being 'leaky' or imperfect, were examined; that is, the immunity provided by the vaccine is not perfect, and vaccinated individuals may still become infected, albeit with a lower probability compared to unvaccinated individuals. This was addressed by incorporating parameters for subgroups with different behaviors or characteristics related to the disease.

The difference of this study from other studies is that it examines the effects of a defective vaccine and mutation on the spread of the COVID-19 outbreak using data obtained for Turkey.

With this study:

(i) Recognizing that a single mathematical model cannot fully represent a real event or concept, we developed a new model by incorporating quantities parameters (mutation, leaky vaccine), new subpopulations (the class infected with the mutant virus) that affect our model and whose behavior is considered for study, differentiating it from existing literature.

(ii) By defining subpopulations (the class infected with the original virus, the class infected with the mutant virus) beyond the compartments in traditional infectious disease models, we addressed the impact of the disease's future progression on these compartments.

(iii) We analyzed the simulations of our model to forecast how interventions, such as vaccination, might influence the course of the disease in efforts to control it.

We organized our study as follows. In the second section, we formulated a mathematical model for the transmission dynamics of COVID-19 infection. We obtained the disease-free equilibrium point and the endemic equilibrium point, and calculated the basic reproduction number. In the third section, we addressed the local stability analysis of these equilibrium points. In the fourth section, we examined the global stability of the disease-free and endemic equilibrium points using the Lyapunov function. In the fifth section, we performed a backward bifurcation analysis of the model using the Castillo-Chavez and Song theorem. In the sixth section, we estimated the parameters of the model. In the seventh section, we conducted a sensitivity analysis of the model. The eighth section encompasses various numerical simulations. In the ninth section, the results are discussed.

## 2 Mathematical model

In this study, we developed a deterministic compartmental model to examine the dynamics of COVID-19 infection. Given that the total population is $N$, the model consists of the following compartments: $S$ (Susceptible), $V$ (Vaccinated), $I_1$ (Infected with Original Virus), $I_2$ (Infected with Mutant Virus), and $R$ (Recovered).

We present below the nonlinear ordinary differential equation system developed to examine the dynamics of COVID-19 infection, along with explanations of the parameters involved in this system.

$$\frac{dS}{dt} = B - \omega S - \beta_1 S I_1 - \beta_2 S I_2 - \alpha S + \mu V + \delta R$$

$$\frac{dV}{dt} = -\omega V - (1-\sigma)\beta_1 V I_1 - (1-\sigma)\beta_2 V I_2 + \alpha S - \mu V$$

$$\frac{dI_1}{dt} = -\omega I_1 - \omega_1 I_1 - m I_1 - r_1 I_1 + \beta_1 S I_1 + (1-\sigma)\beta_1 V I_1 \qquad (2.1)$$

$$\frac{dI_2}{dt} = -\omega I_2 - \omega_2 I_2 + m I_1 - r_2 I_2 + \beta_2 S I_2 + (1-\sigma)\beta_2 V I_2$$

$$\frac{dR}{dt} = -\omega R + r_1 I_1 + r_2 I_2 - \delta R$$

In the model developed in this study; The susceptible population ($S$) increases with birth, loss of vaccine-induced immunity, and loss of natural immunity acquired through infection. This population declines through vaccination (transition to the $V$ class), infection, and natural death. The vaccinated population ($V$) increases by vaccinating susceptible individuals. Since the defective vaccine was present in the model, we took into account that vaccinated people could also become infected. Thus, the vaccinated population decreases due to infection, decrease of artificial active immunity and natural death. The infected population declines with natural death and the population recovering from the disease. The recovered population increases as individuals recover from infection and decreases as immunity acquired by infection declines.

**Table 2.1** Biological interpretations of the parameters

| Parameter | Description |
|---|---|
| $B$ | Natural birth rate |
| $\omega$ | Natural death rate |
| $\beta_1$ | Transmission rate of the original virus |
| $\beta_2$ | Transmission rate of the mutant virus |
| $\alpha$ | Vaccination rate of susceptible individuals |
| $\mu$ | Rate of decline in vaccine-induced immunity |
| $\delta$ | Rate of decline in immunity gained through infection |
| $\sigma$ | Vaccine efficacy |
| $\omega_1$ | Additional mortality rate caused by the original virus |
| $\omega_2$ | Additional mortality rate caused by the mutant virus |
| $m$ | Mutation rate |
| $r_1$ | Recovery rate for the original virus |
| $r_2$ | Recovery rate for the mutant virus |

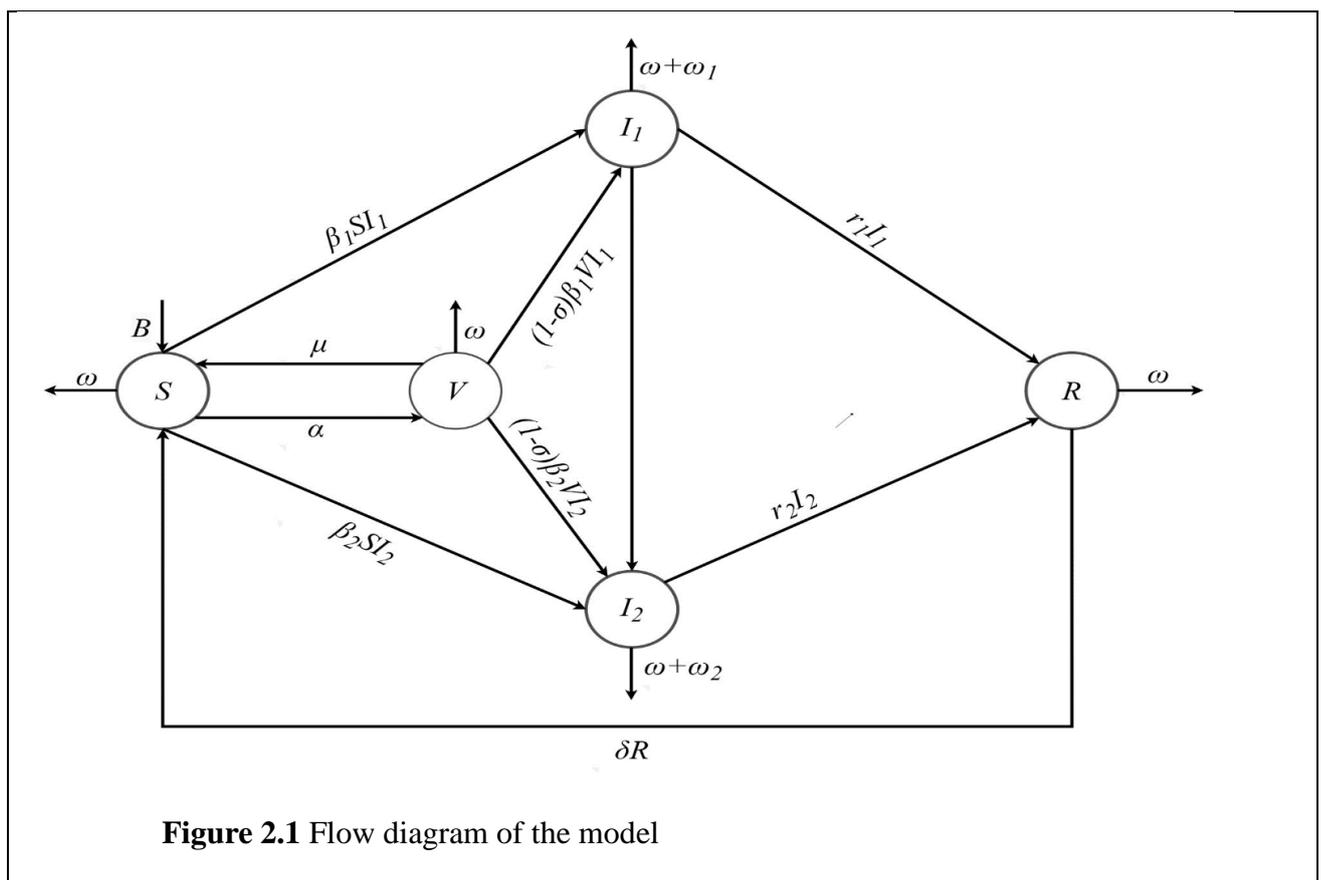

**Figure 2.1** Flow diagram of the model

Initial conditions of the differential equation system for the model;

$$S(0) \geq 0, V(0) \geq 0, I_1(0) \geq 0, I_2(0) \geq 0, R(0) \geq 0$$

is taken as. By summing all the equations in the system,

$$N' = S' + V' + I_1' + I_2' + R'$$

from the equality,

$$N' = B - \omega(S + V + I_1 + I_2 + R) - \omega_1 I_1 - \omega_2 I_2 = B - \omega N - \omega_1 I_1 - \omega_2 I_2 \leq B - \omega N$$

is obtained. Thus,

$$D = \left\{(S, V, I_1, I_2, R) \in \mathbb{R}_+^6 : N \leq \frac{B}{\omega}\right\}$$

which is referred to as the biologically defined region of the mathematical model.

## 2.1 Disease-free equilibrium and endemic equilibrium

The disease-free equilibrium of the system (2.1),

$$E_0 = (S_0, V_0, 0, 0, 0)$$

this is represented as. The disease-free equilibrium

$$E_0 = (S_0, V_0, 0, 0, 0) = \left(\frac{B(\mu + \omega)}{\omega(\mu + \omega + \alpha)}, \frac{\alpha B}{\omega(\mu + \omega + \alpha)}, 0, 0, 0\right)$$

is found as.

Now, let us consider the endemic equilibrium of the system (2.1),

$$E^e = (S^*, V^*, I_1^*, I_2^*, R^*)$$

we can represent it as. The endemic equilibrium,

$$E^e = (S^*, V^*, I_1^*, I_2^*, R^*) = \left(\frac{[B(\omega+\delta)+\delta(r_1 I_1^* + r_2 I_2^*)]\, C}{(\omega+\delta)[(\omega+\beta_1 I_1^* + \beta_2 I_2^* + \alpha)\, C - \mu\alpha]}, \frac{[B(\omega+\delta)+\delta(r_1 I_1^* + r_2 I_2^*)]\, \alpha}{(\omega+\delta)[(\omega+\beta_1 I_1 + \beta_2 I_2 + \alpha)\, C - \mu\alpha]}, I_1^*, I_2^*, \frac{r_1 I_1^* + r_2 I_2^*}{\omega+\delta}\right)$$

is found as. Here,

$$C = \omega + (1-\sigma)\beta_1 I_1^* + (1-\sigma)\beta_2 I_2^* + \mu.$$

## 2.2 Basic reproduction number

In the second edition of his book *The Prevention of Malaria* published in 1911, Ross developed mathematical models of malaria transmission and derived a threshold value, now known as the basic reproduction number [14]. The basic reproduction number, denoted by $\mathcal{R}_0$, is defined as the number of secondary cases generated by a single infected individual in a completely susceptible population. To calculate the basic reproduction number in a compartmental model, a matrix called the next-generation matrix is constructed, which was introduced by Diekmann and Heesterbeek in 1990 [11]. In this study, we used the approach proposed by Van den Driessche and Watmough to construct the next-generation matrix [21].

Let $F$ denote the Jacobian matrix of the terms representing the rate of new infections in the infected compartments at the disease-free equilibrium, and let $V$ denote the Jacobian matrix of the remaining transition terms at the disease-free equilibrium. In this case, for the model (2.1), the matrices $F$ and $V$ are given by:

$$F = \begin{bmatrix} \beta_1 S_0 + (1-\sigma)\beta_1 V_0 & 0 \\ 0 & \beta_2 S_0 + (1-\sigma)\beta_2 V_0 \end{bmatrix}$$

and

$$V = \begin{bmatrix} \omega + \omega_1 + m + r_1 & 0 \\ -m & \omega + \omega_2 + r_2 \end{bmatrix}$$

is obtained. For the system (2.1), the next-generation matrix is

$$K = FV^{-1} = \begin{bmatrix} \dfrac{\beta_1 S_0 + (1-\sigma)\beta_1 V_0}{\omega + \omega_1 + m + r_1} & 0 \\ \dfrac{m(\beta_2 S_0 + (1-\sigma)\beta_2 V_0)}{(\omega + \omega_1 + m + r_1)(\omega + \omega_2 + r_2)} & \dfrac{\beta_2 S_0 + (1-\sigma)\beta_2 V_0}{(\omega + \omega_2 + r_2)} \end{bmatrix}$$

in this form. The eigenvalues of this matrix are

$$\lambda_1 = \frac{\beta_1 S_0 + (1-\sigma)\beta_1 V_0}{\omega + \omega_1 + m + r_1} = \frac{\beta_1 B(\mu + \omega) + (1-\sigma)\beta_1 \alpha B}{(\omega + \omega_1 + m + r_1)\omega(\mu + \omega + \alpha)}$$

$$\lambda_2 = \frac{\beta_2 S_0 + (1-\sigma)\beta_2 V_0}{(\omega + \omega_2 + r_2)} = \frac{\beta_2 B(\mu + \omega) + (1-\sigma)\beta_2 \alpha B}{(\omega + \omega_2 + r_2)\omega(\mu + \omega + \alpha)}$$

$\mathcal{R}_{01}$ represents the basic reproduction number for individuals infected with the original virus, and $\mathcal{R}_{02}$ represents the basic reproduction number for individuals infected with the mutant virus. We obtained is

$$\mathcal{R}_{01} = \frac{\beta_1 B(\mu + \omega) + (1-\sigma)\beta_1 \alpha B}{(\omega + \omega_1 + m + r_1)\omega(\mu + \omega + \alpha)}$$

$$\mathcal{R}_{02} = \frac{\beta_2 B(\mu + \omega) + (1-\sigma)\beta_2 \alpha B}{(\omega + \omega_2 + r_2)\omega(\mu + \omega + \alpha)}$$

## 3. Local stability analysis

In this section, we examined the local stability analysis of the disease-free equilibrium and endemic equilibrium points of the differential equation system representing the mathematical model.

### 3.1 Local asymptotic stability of the disease-free equilibrium

**Theorem 3.1** The disease-free equilibrium point of Model (2.1) is locally asymptotically stable when $\mathcal{R}_{01} < 1$ and $\mathcal{R}_{02} < 1$.

*Proof* The local stability analysis of a mathematical model is conducted using the Jacobian matrix. The sign of the eigenvalues, which are the roots of the characteristic polynomial formed by the Jacobian matrix, provides information about the stability of the nonlinear dynamic system. Considering the Jacobian matrix of the mathematical model (2.1) developed for COVID-19 at the disease-free equilibrium point $E_0$:

$$J_0 = \begin{bmatrix} f_{11} & f_{12} & f_{13} & f_{14} & f_{15} \\ f_{21} & f_{22} & f_{23} & f_{24} & 0 \\ 0 & 0 & f_{33} & 0 & 0 \\ 0 & 0 & f_{43} & f_{44} & 0 \\ 0 & 0 & f_{53} & f_{54} & f_{55} \end{bmatrix} \qquad (3.1)$$

where,

$f_{11} = -\omega - \alpha$, $f_{12} = \mu$, $f_{13} = -\frac{\beta_1 B (\mu+\omega)}{\omega(\mu+\omega+\alpha)}$, $f_{14} = -\frac{\beta_2 B(\mu+\omega)}{\omega(\mu+\omega+\alpha)}$, $f_{15} = \delta$, $f_{21} = \alpha$, $f_{22} = -\omega - \mu$, $f_{23} = -\frac{(1-\sigma)\beta_1 \alpha B}{\omega(\mu+\omega+\alpha)}$, $f_{24} = -\frac{(1-\sigma)\beta_2 \alpha B}{\omega(\mu+\omega+\alpha)}$, $f_{25} = 0$, $f_{31} = 0$, $f_{32} = 0$, $f_{33} = -\omega - \omega_1 - m - r_1 + \frac{\beta_1(B(\mu+\omega)+\alpha B(1-\sigma))}{\omega(\mu+\omega+\alpha)}$, $f_{34} = 0$, $f_{35} = 0$, $f_{43} = m$, $f_{44} = -\omega - \omega_2 - r_2 + \frac{\beta_2(B(\mu+\omega)+\alpha B(1-\sigma))}{\omega(\mu+\omega+\alpha)}$, $f_{53} = r_1$, $f_{54} = r_2$, $f_{55} = -\omega - \delta$.

The eigenvalues, which are the roots of the characteristic polynomial of this matrix are

$\tau_1 = -\omega$, $\tau_2 = -(\omega + \alpha + \mu)$, $\tau_3 = (1 - \mathcal{R}_{01})(-\omega - \omega_1 - m - r_1)$ $\tau_4 = (1 - \mathcal{R}_{02})(-\omega - \omega_2 - r_2)$ ve $\tau_5 = -\omega - \delta$.

The eigenvalues $\tau_{1,2,5}$ have negative real parts. Additionally, under the condition

$$\mathcal{R}_{01} < 1 \,, \mathcal{R}_{02} < 1$$

the eigenvalues $\tau_{3,4}$ also have negative real parts. Thus the disease-free equilibrium point is locally asymptotically stable.

## 3.2 Local asymptotic stability of the endemic equilibrium

**Theorem 3.2** The endemic equilibrium point of Model (2.1) is locally asymptotically stable when $\mathcal{R}_{01} > 1$ and $\mathcal{R}_{02} > 1$.

*Proof* Considering the Jacobian matrix of the mathematical model (2.1) developed for COVID-19 at the endemic equilibrium point $E^e$:

$$J_e = \begin{bmatrix} d_{11} & d_{12} & d_{13} & d_{14} & d_{15} \\ d_{21} & d_{22} & d_{23} & d_{24} & 0 \\ d_{31} & d_{32} & d_{33} & 0 & 0 \\ d_{41} & d_{42} & d_{43} & d_{44} & 0 \\ 0 & 0 & d_{53} & d_{54} & d_{55} \end{bmatrix} \qquad (3.2)$$

where,

$d_{11} = -\omega - \beta_1 I_1^* - \beta_2 I_2^* - \alpha$, $d_{12} = \mu$,

$d_{13} = -\beta_1 S^* = -\beta_1 \frac{[B(\omega+\delta)+\delta(r_1 I_1^* + r_2 I_2^*)]\, C}{(\omega+\delta)[(\omega+\beta_1 I_1^* + \beta_2 I_2^* + \alpha)\, C - \mu\alpha]}$

$d_{14} = -\beta_2 S^* = -\beta_2 \frac{[B(\omega+\delta)+\delta(r_1 I_1^* + r_2 I_2^*)]\, C}{(\omega+\delta)[(\omega+\beta_1 I_1^* + \beta_2 I_2^* + \alpha)\, C - \mu\alpha]}$,

$d_{15} = \delta$, $d_{21} = \alpha$, $d_{22} = -\omega - (1-\sigma)\beta_1 I_1^* - (1-\sigma)\beta_2 I_2^* - \mu$,

$d_{23} = -(1-\sigma)\beta_1 V^* = -(1-\sigma)\beta_1 \frac{[B(\omega+\delta)+\delta(r_1 I_1^* + r_2 I_2^*)]\, \alpha}{(\omega+\delta)[(\omega+\beta_1 I_1^* + \beta_2 I_2^* + \alpha)\, C - \mu\alpha]}$

$d_{24} = -(1-\sigma)\beta_2 V^* = -(1-\sigma)\beta_2 \frac{[B(\omega+\delta)+\delta(r_1 I_1^* + r_2 I_2^*)]\, \alpha}{(\omega+\delta)[(\omega+\beta_1 I_1^* + \beta_2 I_2^* + \alpha)\, C - \mu\alpha]}$, $d_{25} = 0$,

$d_{31} = \beta_1 I_1^*$, $d_{32} = (1-\sigma)\beta_1 I_1^*$, $d_{33} - \omega - \omega_1 - m - r_1 + \beta_1 \frac{[B(\omega+\delta)+\delta(r_1 I_1^* + r_2 I_2^*)]}{(\omega+\delta)[(\omega+\beta_1 I_1 + \beta_2 I_2 + \alpha)\, C - \mu\alpha]}[C + (1-\sigma)\alpha]$, $d_{34} = 0$, $d_{35} = 0$, $d_{41} = \beta_2 I_2^*$, $d_{42} = (1-\sigma)\beta_2 I_2^*$, $d_{43} = m$, $d_{44} = -\omega - \omega_2 - r_2 + \beta_2 \frac{[B(\omega+\delta)+\delta(r_1 I_1^* + r_2 I_2^*)]}{(\omega+\delta)[(\omega+\beta_1 I_1 + \beta_2 I_2 + \alpha)\, C - \mu\alpha]}[C + (1-\sigma)\alpha]$, $d_{45} = 0$, $d_{51} = 0$, $d_{52} = 0$, $d_{53} = r_1$, $d_{54} = r_2$, $d_{55} = -\omega - \delta$.

We can then construct the characteristic polynomial using this matrix:

$$\det(J_e - \tau\, I) = \tau^5 + k_1 \tau^4 + k_2 \tau^3 + k_3 \tau^2 + k_4 \tau + k_5 \qquad (3.3)$$

Here, $\tau$ denotes the eigenvalues, and $k_i$, $i = 1, \dots 5$, represents the coefficients. The coefficients $k_i$ are provided in Appendix A. Thus, if the polynomial (3.3) satisfies the Routh-Hurwitz criteria given by:

$k_1 > 0$, $k_2 > 0$, $k_3 > 0$, $k_4 > 0$, $k_5 > 0$

$$k_1 k_2 k_3 > k_3^2 + k_1^2 k_4$$

$$(k_1 k_4 - k_5)(k_1 k_2 k_3 - k_3^2 - k_1^2 k_4) > k_5 (k_1 k_2 - k_3)^2 + k_1 k_5^2$$

then the endemic equilibrium point is locally asymptotically stable.

**4. Global stability**

In this section, we first examined the global asymptotic stability of the disease-free equilibrium of the model (2.1) using the Lyapunov-Krasovskii-LaSalle Stability Theorem. Subsequently, we investigated the global stability of the endemic equilibrium using the same theorem.

**4.1. Global asymptotic stability of the disease-free equilibrium**

**Theorem 4.1** The disease-free equilibrium of model (2.1) is globally asymptotically stable when $\mathcal{R}_{01} < 1$ and $\mathcal{R}_{02} < 1$.

*Proof* To show that the disease-free equilibrium is globally asymptotically stable, we first define a linear Lyapunov function as follows:

$$L_0(t, S, V, I_1, I_2, R) = C_1 I_1 + C_2 I_2.$$

By taking the derivative of this function and substituting the expressions for $I_1'$ and $I_2'$ from the system (2.1), we obtain

$$\frac{dL_0}{dt} = C_1 I_1' + C_2 I_2'$$

$$= C_1(-\omega I_1 - \omega_1 I_1 - m I_1 - r_1 I_1 + \beta_1 S I_1 + (1-\sigma)\beta_1 V I_1)$$

$$+ C_2(-\omega I_2 - \omega_2 I_2 + m I_1 - r_2 I_2 + \beta_2 S I_2 + (1-\sigma)\beta_2 V I_2)$$

$$= (\mathcal{R}_{01} - 1)I_1 + (\mathcal{R}_{02} - 1)I_2 \leq 0$$

where,

$$C_1 = \frac{1}{\omega + \omega_1 + m + r_1} - \frac{m\omega(\mu + \omega + \alpha)}{(\omega + \omega_2 + r_2)[-(\omega + \omega_1 + m + r_1)\omega(\mu + \omega + \alpha) + \beta_1 B(\mu + \omega) + (1-\sigma)\beta_1 \alpha B]},$$

$$C_2 = \frac{1}{\omega + \omega_2 + r_2}.$$

Thus, according to the Lyapunov Theorem, the disease-free equilibrium $E_0$ is globally asymptotically stable when $\mathcal{R}_{01} < 1$ and $\mathcal{R}_{02} < 1$.

**4.2 Global asymptotic stability of the endemic equilibrium**

**Theorem 4.2** The endemic equilibrium of model (2.1) is globally asymptotically stable when $\mathcal{R}_{01} > 1$ and $\mathcal{R}_{02} > 1$.

*Proof* To determine the global stability of the endemic equilibrium, we introduce the following Lyapunov function:

$$L(t, S, V, I_1, I_2, R) = \kappa_1 \left(S - S^* - S^* \ln \frac{S}{S^*}\right) + \kappa_2 \left(V - V^* - V^* \ln \frac{V}{V^*}\right) + \kappa_3 \left(I_1 - I_1^* - I_1^* \ln \frac{I_1}{I_1^*}\right) + \kappa_4 \left(I_2 - I_2^* - I_2^* \ln \frac{I_2}{I_2^*}\right) + \kappa_5 \left(R - R^* - R^* \ln \frac{R}{R^*}\right).$$

Here, $\kappa_1 > 0$, $\kappa_2 > 0$, $\kappa_3 > 0$, $\kappa_4 > 0$ and $\kappa_5 > 0$. When $(S, V, I_1, I_2, R) = (S^*, V^*, I_1^*, I_2^*, R^*)$, $L = 0$, and otherwise $L > 0$. Additionally, $L$ is radially unbounded and

$$L'(x) = \frac{d}{dt} L(x(t)) = \frac{\partial L}{\partial x} \frac{dx}{dt}.$$

Let's show that the derivative of $L$ with respect to $t$ is negative. The derivative of the Lyapunov function is given by

$$\frac{dL}{dt} = \kappa_1 \left(1 - \frac{S^*}{S}\right) S' + \kappa_2 \left(1 - \frac{V^*}{V}\right) V' + \kappa_3 \left(1 - \frac{I_1^*}{I_1}\right) I_1' + \kappa_4 \left(1 - \frac{I_2^*}{I_2}\right) I_2' + \kappa_5 \left(1 - \frac{R^*}{R}\right) R'.$$

Thus, by substituting the derivatives $S', V', I_1', I_2'$, and $R'$ from the system (2.1) and performing some simplifications, we obtain

$$\frac{dL}{dt} = -(\omega + \alpha)\kappa_1\left(1 - \frac{1}{x_1}\right)^2 S + \kappa_1\left(1 - \frac{1}{x_1} - x_3 x_1 + x_3\right)\beta_1 S^* I_1^*$$

$$+ \kappa_1\left(1 - \frac{1}{x_1} - x_4 x_1 + x_4\right)\beta_2 S^* I_2^* + \kappa_1\left(x_5 - \frac{x_5}{x_1} - 1 + \frac{1}{x_1}\right)\delta R^*$$

$$+ \kappa_1\left(x_2 - \frac{x_2}{x_1} - 1 + \frac{1}{x_1}\right)\mu V^* + \left(-x_2 x_3 + x_3 + 1 - \frac{1}{x_2}\right)\beta_1 V^* I_1^*(1-\sigma)\kappa_2$$

$$+ \left(-x_2 x_4 + x_4 + 1 - \frac{1}{x_2}\right)\beta_2 V^* I_2^*(1-\sigma)\kappa_2$$

$$+ \left(-1 + \frac{1}{x_3} + x_1 x_3 - x_1\right)\beta_1 S^* I_1^* \kappa_3$$

$$+ \left(-1 + \frac{1}{x_3} + x_2 x_3 - x_2\right)\beta_1 V^* I_1^*(1-\sigma)\kappa_3$$

$$+ \left(-1 + \frac{1}{x_4} + x_1 x_4 - x_1\right)\beta_2 S^* I_2^* \kappa_4$$

$$+ \left(-1 + \frac{1}{x_4} + x_2 x_4 - x_2\right)\beta_2 V^* I_2^*(1-\sigma)\kappa_4 + \left(-1 + \frac{1}{x_5} + x_3 - \frac{x_3}{x_5}\right)r_1 I_1^* \kappa_5$$

$$+ \left(-1 + \frac{1}{x_5} + x_4 - \frac{x_4}{x_5}\right)r_2 I_2^* \kappa_5,$$

where,

$$x_1 = \frac{S}{S^*}, x_2 = \frac{V}{V^*}, x_3 = \frac{I_1}{I_1^*}, x_4 = \frac{I_2}{I_2^*}, x_5 = \frac{R}{R^*}$$

is obtained, and

$$
\begin{aligned}
T = {} & \kappa_1\beta_1 S^* I_1^* + \kappa_1\beta_2 S^* I_2^* - \kappa_1\delta R^* - \kappa_1\mu V^* + \beta_1 V^* I_1^*(1-\sigma)\kappa_2 + \beta_2 V^* I_2^*(1-\sigma)\kappa_2 \\
& - \beta_1 S^* I_1^* \kappa_3 - \beta_1 V^* I_1^*(1-\sigma)\kappa_3 - \beta_2 S^* I_2^* \kappa_4 - \beta_2 V^* I_2^*(1-\sigma)\kappa_4 - r_1 I_1^* \kappa_5 \\
& - r_2 I_2^* \kappa_5 + (-\kappa_1\beta_1 S^* I_1^* - \kappa_1\beta_2 S^* I_2^* + \kappa_1\delta R^* + \kappa_1\mu V^*)\frac{1}{x_1} \\
& + (-\beta_1 S^* I_1^* \kappa_3 - \beta_2 S^* I_2^* \kappa_4)x_1 \\
& + (-\beta_1 V^* I_1^*(1-\sigma)\kappa_2 - \beta_2 V^* I_2^*(1-\sigma)\kappa_2)\frac{1}{x_2} \\
& + (\kappa_1\mu V^* - \beta_1 V^* I_1^*(1-\sigma)\kappa_3 - \beta_2 V^* I_2^*(1-\sigma)\kappa_4)x_2 \\
& + (\beta_1 S^* I_1^* \kappa_3 + \beta_1 V^* I_1^*(1-\sigma)\kappa_3)\frac{1}{x_3} \\
& + (\kappa_1\beta_1 S^* I_1^* + \beta_1 V^* I_1^*(1-\sigma)\kappa_2 + r_1 I_1^* \kappa_5)x_3 \\
& + (\beta_2 S^* I_2^* \kappa_4 + \beta_2 V^* I_2^*(1-\sigma)\kappa_4)\frac{1}{x_4} \\
& + (\kappa_1\beta_2 S^* I_2^* + \beta_2 V^* I_2^*(1-\sigma)\kappa_2 + r_2 I_2^* \kappa_5)x_4 + (r_1 I_1^* \kappa_5 + r_2 I_2^* \kappa_5)\frac{1}{x_5} \\
& + (\kappa_1\delta R^*)x_5 + (-\kappa_1\beta_2 S^* I_2^* + \beta_2 S^* I_2^* \kappa_4)x_4 x_1 \\
& + (-\beta_1 V^* I_1^*(1-\sigma)\kappa_2 + \beta_1 V^* I_1^*(1-\sigma)\kappa_3)x_2 x_3 \\
& + (-\beta_2 V^* I_2^*(1-\sigma)\kappa_2 + \beta_2 V^* I_2^*(1-\sigma)\kappa_4)x_2 x_4 \\
& + (-\kappa_1\beta_1 S^* I_1^* + \beta_1 S^* I_1^* \kappa_3)x_1 x_3 - \kappa_1\delta R^* \frac{x_5}{x_1} - \kappa_1\mu V^* \frac{x_2}{x_1} - r_1 I_1^* \kappa_5 \frac{x_3}{x_5} \\
& - r_2 I_2^* \kappa_5 \frac{x_4}{x_5}
\end{aligned}
$$

considering that we have

$$
\frac{dL}{dt} = -(\omega + \alpha)\kappa_1 \left(1 - \frac{1}{x_1}\right)^2 S + T(x_1, x_2, x_3, x_4, x_5)
$$

where,

$$
\kappa_1 = \kappa_2 = \kappa_3 = \kappa_4 = 1
$$

this implies that all the coefficients $x_4 x_1$, $x_2 x_3$, $x_2 x_4$, $x_1 x_3$ are zero. Let's rewrite $T$ as follows:

$$T = -\left(-\frac{1}{x_3} + x_1 + \frac{1}{x_1} - x_3\right)\beta_1 S^* I_1^* - \left(-\frac{1}{x_4} + x_1 + \frac{1}{x_1} - x_4\right)\beta_2 S^* I_2^*$$

$$- \left(-x_5 + \frac{x_5}{x_1} + 1 - \frac{1}{x_1}\right)\delta R^* - \left(-x_2 + \frac{x_2}{x_1} + 1 - \frac{1}{x_1}\right)\mu V^*$$

$$- \left(-\frac{1}{x_3} + x_2 - x_3 + \frac{1}{x_2}\right)\beta_1 V^* I_1^* (1-\sigma)$$

$$- \left(-\frac{1}{x_4} + x_2 - x_4 + \frac{1}{x_2}\right)\beta_2 V^* I_2^* (1-\sigma) - \left(1 - \frac{1}{x_5} - x_3 + \frac{x_3}{x_5}\right)r_1 I_1^* \kappa_5$$

$$- \left(1 - \frac{1}{x_5} - x_4 + \frac{x_4}{x_5}\right)r_2 I_2^* \kappa_5.$$

If we take $x_1 = 1$, $x_2 = 1$, $x_3 = 1$, $x_4 = 1$ and $x_5 = 1$, then $T \leq 0$. Thus, we obtain

$$\frac{dL}{dt} = -(\omega + \alpha)\kappa_1\left(1 - \frac{1}{x_1}\right)^2 S + T(x_1, x_2, x_3, x_4, x_5) \leq 0$$

According to LaSalle's Theorem, the invariant set is defined by

$$\Omega = \{x \in \mathbb{R}^n : L'(x) = 0\}$$

Since the invariant set contains only the endemic equilibrium $(S^*, V^*, I_1^*, I_2^*, R^*)$, we conclude that the endemic equilibrium is globally asymptotically stable.

## 5. Backward bifurcation

The qualitative dynamics of epidemiological models depend on the threshold value known as the basic reproduction number, denoted as $\mathcal{R}_0$. Backward bifurcation is examined by varying the bifurcation parameter, which is identified as the basic reproduction number. It should be noted that $\mathcal{R}_0$ is a sum of parameters in the model. In general, if $\mathcal{R}_0 < 1$, the disease will eventually disappear. If $\mathcal{R}_0 > 1$, an infected individual can transmit the disease to more than one person, allowing the disease to continue spreading within the population. However, in some cases, despite $\mathcal{R}_0 < 1$, the disease may persist in transmission. In such scenarios, the phenomenon of backward bifurcation arises [21]. In this study, we conducted a bifurcation analysis using the Castillo-Chavez and Song Bifurcation Theorem.

Initially, we demonstrated the presence of backward bifurcation for the original virus. For this purpose, we selected $\beta_1$ as the bifurcation parameter and determined $\beta_1$ by setting $R_{01} = 1$. In this case

$$\beta_1 = \beta_1^* = \frac{\omega(\omega + \omega_1 + m + r_1)(\mu + \omega + \alpha)}{B(\mu + \omega) + (1 - \sigma)\alpha B}$$

we obtain the value of this parameter. By substituting this parameter value into the Jacobian matrix $J_0$, given in equation (3.1), we obtain

$$J_{(E_0, \beta_1^*)} = \begin{bmatrix} g_{11} & g_{12} & g_{13} & g_{14} & g_{15} \\ g_{21} & g_{22} & g_{23} & g_{24} & 0 \\ 0 & 0 & 0 & 0 & 0 \\ 0 & 0 & g_{43} & g_{44} & 0 \\ 0 & 0 & g_{53} & g_{54} & g_{55} \end{bmatrix}$$

Here,

$g_{11} = -\omega - \alpha,\ g_{12} = \mu,\ g_{13} = -\frac{(\omega + \omega_1 + m + r_1)(\mu + \omega)}{(\mu + \omega) + (1 - \sigma)\alpha},\ g_{14} = -\frac{\beta_2 B(\mu + \omega)}{\omega(\mu + \omega + \alpha)},\ g_{15} = \delta,$

$g_{21} = \alpha,\ g_{22} = -\omega - \mu,\ g_{23} = -\frac{(\omega + \omega_1 + m + r_1)(1 - \sigma)\alpha}{(\mu + \omega) + (1 - \sigma)\alpha},\ g_{24} = -\frac{(1 - \sigma)\beta_2 \alpha B}{\omega(\mu + \omega + \alpha)},\ g_{33} = 0,\ g_{43} = m,\ g_{44} = -\omega - \omega_2 - r_2 + \frac{\beta_2(B(\mu + \omega) + \alpha B(1 - \sigma))}{\omega(\mu + \omega + \alpha)},\ g_{53} = r_1,\ g_{54} = r_2,\ g_{55} = -\omega - \delta.$

The characteristic polynomial and eigenvalues of the matrix $J_{(E_0, \beta_1^*)}$ are as follows:

$$|J_{(E_0, \beta_1^*)} - \tau I| = [(g_{11} - \tau)(g_{22} - \tau) - g_{21} g_{12}][(-\tau)(g_{44} - \tau)(g_{55} - \tau)]$$

$$= [\tau^2 + (2\omega + \alpha + \mu)\tau + (\omega^2 + \omega\mu + \alpha\omega)](-\tau)\left(-\omega - \omega_2 - r_2 + \frac{\beta_2(B(\mu + \omega) + \alpha B(1 - \sigma))}{\omega(\mu + \omega + \alpha)} - \tau\right)(-\omega - \delta - \tau),$$

$\tau_1 = -\omega,\ \tau_2 = -(\omega + \alpha + \mu),\ \tau_3 = 0,\ \tau_4 = -\omega - \omega_2 - r_2 + \frac{\beta_2(B(\mu + \omega) + \alpha B(1 - \sigma))}{\omega(\mu + \omega + \alpha)},\ \tau_5 = -\omega - \delta.$

The matrix $J_{(E_0, \beta_1^*)}$ has a simple zero eigenvalue, while the remaining eigenvalues are negative. Now, let us compute the right eigenvector $\mathbf{w} = (w_1, w_2, w_3, w_4, w_5)$ corresponding to the zero eigenvalue, as well as the left eigenvector $\mathbf{v} = (v_1, v_2, v_3, v_4, v_5)$.

$$J_{(E_0, \beta_1^*)} \mathbf{w} = 0.$$

By solving the equation, we obtain the components of the right eigenvector $\mathbf{w} = (w_1, w_2, w_3, w_4, w_5)$,

$$w_1 = \frac{(\omega + \mu)X + \mu(\omega + \delta)Y}{m(\omega + \delta)[(\mu + \omega) + (1 - \sigma)\alpha][-\omega(\omega + \mu + \alpha)]},$$

$$w_2 = \frac{\alpha X + (\omega + \alpha)(\omega + \delta)Y}{m(\omega + \delta)[(\mu + \omega) + (1 - \sigma)\alpha][-\omega(\omega + \mu + \alpha)]},$$

$$w_3 = \frac{(1 - R_{02})}{m} A_1,$$

$$w_4 = \omega(\mu + \omega + \alpha),$$

$$w_5 = \frac{r_1(1 - R_{02})A_1 + r_2 m \omega(\mu + \omega + \alpha)}{m(\omega + \delta)},$$

where,

$A_1 = \omega(\mu + \omega + \alpha)(\omega + \omega_2 + r_2),$

$X = (\omega + \omega_1 + m + r_1)(\mu + \omega)(\omega + \delta)(1 - R_{02})A_1$
$\quad + \beta_2 B(\mu + \omega)m(\omega + \delta)[(\mu + \omega) + (1 - \sigma)\alpha]$
$\quad - \delta[r_1(1 - R_{02})A_1 + r_2 m \omega(\mu + \omega + \alpha)][(\mu + \omega) + (1 - \sigma)\alpha],$

$Y = (\omega + \omega_1 + m + r_1)(1 - \sigma)\alpha(1 - R_{02})A_1 + (1 - \sigma)\beta_2 \alpha B m[(\mu + \omega) + (1 - \sigma)\alpha].$

Similarly, solving the system $\mathbf{v} J_{(E_0, \beta_1^*)} = 0$ yields the components of the left eigenvector $\mathbf{v} = (v_1, v_2, v_3, v_4, v_5)$ as $v_1 = 0$, $v_2 = 0$, $v_3 = 1$, $v_4 = 0$, $v_5 = 0$. Now, let us calculate the bifurcation constants $a$ and $b$ as given in [4]:

$$a = \sum_{k,i,j=1}^{n} v_k w_i w_j \frac{\partial^2 f_k}{\partial x_i \partial x_j}(0,0) \tag{5.1}$$

$$b = \sum_{k,i=1}^{n} v_k w_i \frac{\partial^2 f_k}{\partial x_i \partial \phi}(0,0) \tag{5.2}$$

Let's rewrite the system (2.1) as follows:

$$\frac{dx_1}{dt} = B - \omega x_1 - \beta_1 x_1 x_3 - \beta_2 x_1 x_4 - \alpha x_1 + \mu x_2 + \delta x_5 = f_1$$

$$\frac{dx_2}{dt} = -\omega x_2 - (1 - \sigma)\beta_1 x_2 x_3 - (1 - \sigma)\beta_2 x_2 x_4 + \alpha x_1 - \mu x_2 = f_2$$

$$\frac{dx_3}{dt} = -\omega x_3 - \omega_1 x_3 - m x_3 - r_1 x_3 + \beta_1 x_1 x_3 + (1 - \sigma)\beta_1 x_2 x_3 = f_3 \tag{5.3}$$

$$\frac{dx_4}{dt} = -\omega x_4 - \omega_2 x_4 + m x_3 - r_2 x_4 + \beta_2 x_1 x_4 + (1-\sigma)\beta_2 x_2 x_4 = f_4$$

$$\frac{dx_5}{dt} = -\omega x_5 + r_1 x_3 + r_2 x_4 - \delta x_5 = f_5.$$

By calculating the required partial derivatives in (5.1) and (5.2) and writing those that are non-zero, we obtain

$$a = -\frac{(1-R_{02})^2 A_1^2 (\omega+\omega_1+m+r_1)^2 (\mu+\omega)}{m^2 B[(\mu+\omega)+(1-\sigma)\alpha]}$$

$$-(1-R_{02})A_1(\omega+\omega_1+m+r_1)(\mu+\omega)\beta_2$$

$$-\frac{[\mu+(1-\sigma)(\omega+\alpha)](\omega+\omega_1+m+r_1)^2(1-\sigma)\alpha(1-R_{02})^2 A_1^2}{m^2 B[(\mu+\omega)+(1-\sigma)\alpha]^2}$$

$$-\frac{[\mu+(1-\sigma)(\omega+\alpha)](1-\sigma)\beta_2\alpha(1-R_{02})A_1(\omega+\omega_1+m+r_1)}{m[(\mu+\omega)+(1-\sigma)\alpha]}$$

$$+\delta\frac{(1-R_{02})A_1(\omega+\omega_1+m+r_1)w_5}{mB}$$

and

$$b = \frac{(1-R_{02})}{m} A_1 [S_0 + (1-\sigma)V_0].$$

Since $R_{02} > 0$, $S_0 > 0$, $V_0 > 0$ and $1-\sigma > 0$, it follows that $b > 0$. To ensure that the condition $a > 0$ is met, the parameter $\delta$ must.

$$\delta > \delta^* = \frac{1}{w_5}\left[\frac{(1-R_{02})A_1(\omega+\omega_1+m+r_1)(\mu+\omega)}{m[(\mu+\omega)+(1-\sigma)\alpha]} + (\mu+\omega)\beta_2 B\right.$$

$$+ \frac{[\mu+(1-\sigma)(\omega+\alpha)](\omega+\omega_1+m+r_1)(1-\sigma)\alpha(1-R_{02})A_1}{m[(\mu+\omega)+(1-\sigma)\alpha]^2}$$

$$\left.+ \frac{[\mu+(1-\sigma)(\omega+\alpha)](1-\sigma)\beta_2\alpha B}{[(\mu+\omega)+(1-\sigma)\alpha]}\right] \quad (5.4)$$

satisfy the condition. Therefore, according to the Castillo-Chavez and Song Bifurcation Theorem, backward bifurcation will occur for $\beta_1$ if the condition given in equation (5.4) is satisfied.

Taking $\mathcal{R}_{02} = 1$ we obtain the bifurcation parameter $\beta_2$ as follows:

$$\beta_2 = \beta_2^* = \frac{\omega(\omega+\omega_2+r_2)(\mu+\omega+\alpha)}{B(\mu+\omega)+(1-\sigma)\alpha B}.$$

By substituting $\beta_2^*$ into the Jacobian matrix $J_0$ at the disease-free equilibrium, we obtain

$$J_{(E_0,\beta_2^*)} = \begin{bmatrix} k_{11} & k_{12} & k_{13} & k_{14} & k_{15} \\ k_{21} & k_{22} & k_{23} & k_{24} & 0 \\ 0 & 0 & k_{33} & 0 & 0 \\ 0 & 0 & k_{43} & 0 & 0 \\ 0 & 0 & k_{53} & k_{54} & k_{55} \end{bmatrix}$$

where,

$k_{11} = -\omega - \alpha$, $k_{12} = \mu$, $k_{13} = -\frac{\beta_1 B\,(\mu+\omega)}{\omega(\mu+\omega+\alpha)}$, $k_{14} = -\frac{(\omega+\omega_2+r_2)(\mu+\omega)}{(\mu+\omega)+(1-\sigma)\alpha}$, $k_{15} = \delta$, $k_{21} = \alpha$,

$k_{22} = -\omega - \mu$, $k_{23} = -\frac{(1-\sigma)\beta_1\alpha B}{\omega(\mu+\omega+\alpha)}$, $k_{24} - \frac{(\omega+\omega_2+r_2)(1-\sigma)\alpha}{(\mu+\omega)+(1-\sigma)\alpha}$, $k_{33} = -\omega - \omega_1 - m - r_1 +$

$\frac{\beta_1(B\,(\mu+\omega)+\alpha B(1-\sigma))}{\omega(\mu+\omega+\alpha)}$, $k_{43} = m$, $k_{44} = 0$, $k_{53} = r_1$, $k_{54} = r_2$, $k_{55} = -\omega - \delta$

is obtain. Let's write the characteristic polynomial and eigenvalues of the matrix $J_{(E_0,\beta_2^*)}$:

$$|J_{(E_0,\beta_2^*)} - \tau I| = \begin{vmatrix} k_{11} - \tau & k_{12} & k_{13} & k_{14} & k_{15} \\ k_{21} & k_{22} - \tau & k_{23} & k_{24} & 0 \\ 0 & 0 & k_{33} - \tau & 0 & 0 \\ 0 & 0 & k_{43} & -\tau & 0 \\ 0 & 0 & k_{53} & k_{54} & k_{55} - \tau \end{vmatrix} = 0$$

$$|J_{(E_0,\beta_2^*)} - \tau I| = [(k_{11} - \tau)(k_{22} - \tau) - k_{21}\,k_{12}]\,[(k_{33} - \tau)(-\tau)(k_{55} - \tau)]$$
$$= [\tau^2 + (2\omega + \alpha + \mu)\tau + (\omega^2 + \omega\mu + \alpha\omega)](k_{33} - \tau)(-\tau)(-\omega - \delta - \tau)$$

and the eigenvalues are:

$\tau_1 = -\omega$, $\tau_2 = -(\omega + \alpha + \mu)$, $\tau_3 = -\omega - \omega_1 - m - r_1 + \frac{\beta_1(B\,(\mu+\omega)+\alpha B(1-\sigma))}{\omega(\mu+\omega+\alpha)}$, $\tau_4 = 0$,

$\tau_5 = -\omega - \delta$

it is observed that the matrix $J_{(E_0,\beta_2^*)}$ has a simple zero eigenvalue at $E_0$ and $\beta_2 = \beta_2^*$, while the remaining eigenvalues are negative.

If the procedures performed for $\beta_1$ are repeated for $\beta_2$,

$w_{33} = 0$, $w_{44} = 1$, $w_{55} = \frac{r_2}{\omega+\delta}$,

$w_{22} = \frac{\alpha(\omega+\omega_2+r_2)(\mu+\omega)(\omega+\delta) - \alpha\delta r_2[(\mu+\omega)+(1-\sigma)\alpha]+(\omega+\alpha)(\omega+\omega_2+r_2)(1-\sigma)\alpha(\omega+\delta)}{[-\omega(\omega+\mu+\alpha)][(\mu+\omega)+(1-\sigma)\alpha](\omega+\delta)}$

$w_{11} = \frac{(\omega+\omega_2+r_2)(\mu+\omega)^2(\omega+\delta) - (\omega+\mu)\delta r_2[(\mu+\omega)+(1-\sigma)\alpha] + \mu(\omega+\omega_2+r_2)(1-\sigma)\alpha(\omega+\delta)}{[-\omega(\omega+\mu+\alpha)][(\mu+\omega)+(1-\sigma)\alpha](\omega+\delta)}$

resulting in the right eigenvector $\boldsymbol{w_1} = (w_{11}, w_{22}, w_{33}, w_{44}, w_{55})$ and

$$v_{11} = 0,\ v_{22} = 0,\ v_{55} = 0,\ v_{44} = 1,\ v_{33} = \frac{-m\omega(\mu+\omega+\alpha)}{-(\omega+\omega_1+m+r_1)\omega(\mu+\omega+\alpha)+\beta_1 B[(\mu+\omega)+\alpha(1-\sigma)]}$$

resulting in the left eigenvector $\boldsymbol{v_1} = (v_{11}, v_{22}, v_{33}, v_{44}, v_{55})$.

In this case, the bifurcation constants $a$ and $b$ are

$$a = -(\omega+\omega_2+r_2)\frac{(\omega+\omega_2+r_2)(\mu+\omega)[(\mu+\omega)+(1-\sigma)\alpha]+(\omega+\omega_2+r_2)(1-\sigma)\alpha[\mu+(\omega+\alpha)(1-\sigma)]}{B[(\mu+\omega)+(1-\sigma)\alpha]^2} + \delta\ w_{55}\frac{(\omega+\omega_2+r_2)}{B}$$

and

$$b = S_0 + (1-\sigma)V_0$$

is obtain. It can be easily seen that $b > 0$.

$$\delta > \delta^* = \frac{1}{w_{55}}\frac{(\omega+\omega_2+r_2)(\mu+\omega)[(\mu+\omega)+(1-\sigma)\alpha]+(\omega+\omega_2+r_2)(1-\sigma)\alpha[\mu+(\omega+\alpha)(1-\sigma)]}{[(\mu+\omega)+(1-\sigma)\alpha]^2} \quad (5.5)$$

given that the constant $a$ will be positive when the condition is satisfied and that the constant $b$ is also positive, backward bifurcation will occur at $\beta_2 = \beta_2^*$ according to the Castillo-Chavez and Song Bifurcation Theorem.

**Theorem 5.1** Consider the system (2.1). A backward bifurcation occurs in system (2.1) if inequalities (5.4) and (5.5) hold at $\mathcal{R}_{01} = 1$ and $\mathcal{R}_{02} = 1$.

## 6 Parameter estimation

Certain parameters, such as birth rates, natural mortality rates, and recovery rates, can be directly inferred from population and epidemiological data. However, other parameters may not be as straightforward to estimate from direct observations. For more complex models or for the estimation of multiple parameters simultaneously, a systematic approach is essential [19]. Consequently, parameter estimation methodologies should be employed to accurately determine these parameters.

In this section, our goal is to estimate the model parameters by fitting the nonlinear ordinary differential equation system we developed to the disease data in order to study the dynamics of COVID-19 infection. In our study, we conducted a model analysis using the data of COVID-19 cases reported in Turkey between December 20, 2021, and January 31, 2022. We solved the system using Python's ***odeint*** function and then estimated the parameters of the $SVI_1I_2R$ model by calculating the error between the model's predictions and the actual data with the help of the *fit_model* function. Using this function, we optimized the model parameters ($\beta_1, \beta_2, m$). Subsequently, we employed the *minimize* function to minimize the error function. Finally, we visualized the solutions in comparison with the data.

In the model we developed for the COVID-19 pandemic, there are 13 different parameter values. We estimated the transmissibility rate of the original virus ($\beta_1$), the transmissibility rate of the mutant virus ($\beta_2$) and the mutation rate ($m$). For the initial values of these parameters, we used $\beta_1 = 0.000000004$, $\beta_2 = 0.0000000085$ and $m = 0.01$. The optimal values for these parameters were found to be $\beta_1 = 0.000000001167817614$, $\beta_2 = 0.000000007368542050$ and $m = 0.009308731535398908$.

**Table 6.1** Initial values of parameters

| Parameters | Description | Initial value | Source |
|---|---|---|---|
| $B$ | Natural birth rate | 2993 | [30] |
| $\omega$ | Natural death rate | 0.00003535 | [31] |
| $\beta_1$ | Transmission rate of the original virus | 0.000000001167817614 | Estimated |
| $\beta_2$ | Transmission rate of the mutant virus | 0.000000007368542050 | Estimated |
| $\alpha$ | Vaccination rate of susceptible individuals | 0.012 | Assumed |
| $\mu$ | Rate of decline in vaccine-induced immunity | 0.005 day | Assumed |
| $\delta$ | Rate of decline in immunity gained through infection | 0.0027 day | [7], [32] |
| $\sigma$ | Vaccine efficacy | 0.9 | Assumed |
| $\omega_1$ | Additional mortality rate caused by the original virus | 0.005579 | [27] |
| $\omega_2$ | Additional mortality rate caused by the mutant virus | 0.002286 | [27] |
| $m$ | Mutation rate | 0.009308731535398908 | Estimated |
| $r_1$ | Recovery rate for the original virus | 0.0833 | [19] |
| $r_2$ | Recovery rate for the mutant virus | 0.1 | [19] |

## 7 Sensitivity analysis

Sensitivity analysis is the study of how uncertainty in a model's output can be attributed to various sources of uncertainty in the model's inputs [28]. In epidemiology, significant static quantities, such as the basic reproductive number, equilibrium prevalence, and equilibrium incidence, are dependent on the parameters of differential equation models. Understanding how these quantities respond to changes in parameters is crucial. The change in the output quantity $Q$ with respect to a parameter $p$ is measured by the derivative of the output with respect to the parameter. Therefore, sensitivity analysis is initially viewed as a local measure of the effect of a specific input on a specific output [21]. This form of sensitivity analysis, known as local sensitivity analysis, is obtained using the partial derivative of the output function with respect to the input factors when the input factors are known with minimal uncertainty [20]. The sensitivity index, which quantifies the change in the output quantity $Q$ with respect to a parameter $p$, is calculated as follows. The sensitivity index is represented by

$$\mathcal{Z}_p^Q = \frac{p}{Q}\frac{\partial Q}{\partial p}$$

In this section, we calculated the sensitivity index for each model parameter related to the basic reproductive number. This index indicates the relative importance of each parameter in the model in terms of its impact on COVID-19 transmission. The index is used to identify the parameter with the greatest effect on the basic reproductive number, which serves as the target for interventions. Parameters that have a significant impact on the basic reproductive number are crucial in managing the COVID-19 outbreak [26]. By calculating the sensitivity indices of $\mathcal{R}_{01}$ and $\mathcal{R}_{02}$ with respect to all model parameters from our mathematical model, using the formulas

$$\mathcal{Z}_p^{\mathcal{R}_{01}} = \frac{p}{\mathcal{R}_{01}}\frac{\partial \mathcal{R}_{01}}{\partial p} \qquad \mathcal{Z}_p^{\mathcal{R}_{02}} = \frac{p}{\mathcal{R}_{02}}\frac{\partial \mathcal{R}_{02}}{\partial p}$$

we obtained the following table.

**Table 7.1** Sensitivity indices of $\mathcal{R}_{01}$ and $\mathcal{R}_{02}$ with respect to all model parameters

| Parameter | Value | Sensitivity index value ($\mathcal{R}_{01}$) | Analysis | Sensitivity index value ($\mathcal{R}_{02}$) | Analysis |
|---|---|---|---|---|---|
| $B$ | 2993 | 1 | Positive | 1 | Positive |
| $\omega$ | 0.00003535 | $-0.9971$ | Negative | $-0.99635$ | Negative |
| $\beta_1$ | 0.000000001167817614 | 1 | Positive | 0 | Neutral |
| $\beta_2$ | 0.000000007368542050 | 0 | Neutral | 1 | Positive |
| $\alpha$ | 0.01795 | $-0.5116$ | Negative | $-0.297$ | Negative |
| $\mu$ | 0.005 | 0.5084 | Positive | 0.5084 | Positive |
| $\delta$ | 0.0027 | 0 | Neutral | 0 | Neutral |
| $\sigma$ | 0.9 | $-1.7320$ | Negative | $-1.7320$ | Negative |
| $\omega_1$ | 0.005579 | $-0.0568$ | Negative | 0 | Neutral |
| $\omega_2$ | 0.002286 | 0 | Neutral | $-0.02234$ | Negative |
| $m$ | 0.009308731535398908 | $-0.0948$ | Negative | 0 | Neutral |
| $r_1$ | 0.0833 | $-0.8481$ | Negative | 0 | Neutral |
| $r_2$ | 0.1 | 0 | Neutral | $-0.9773$ | Negative |

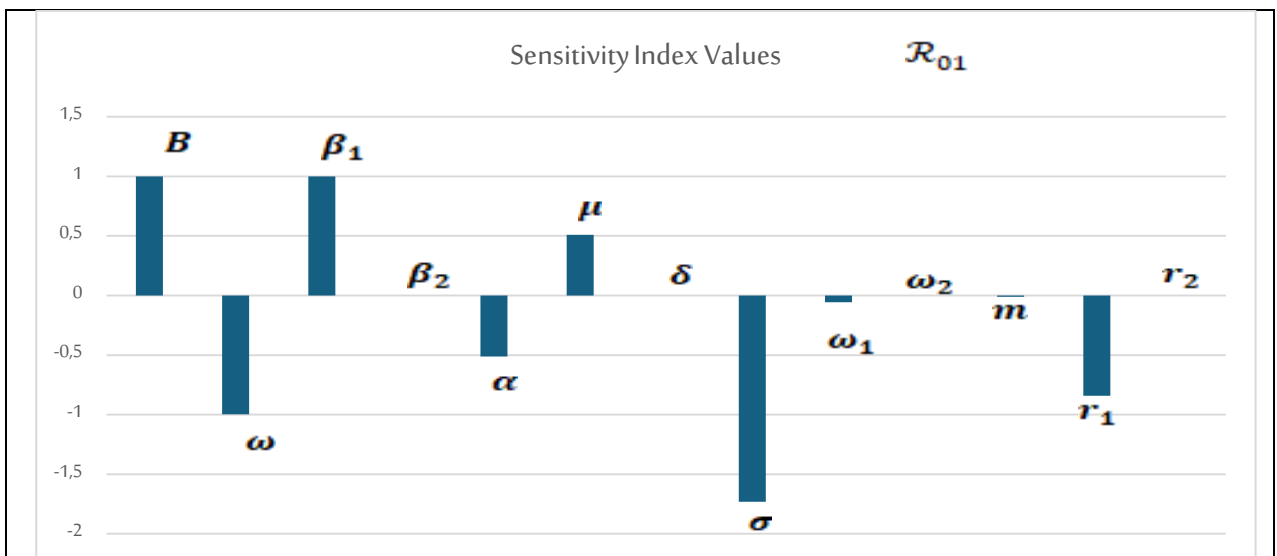

**Figure 7.1** Basic reproduction number $\mathcal{R}_{01}$ and sensitivity index values of the parameters

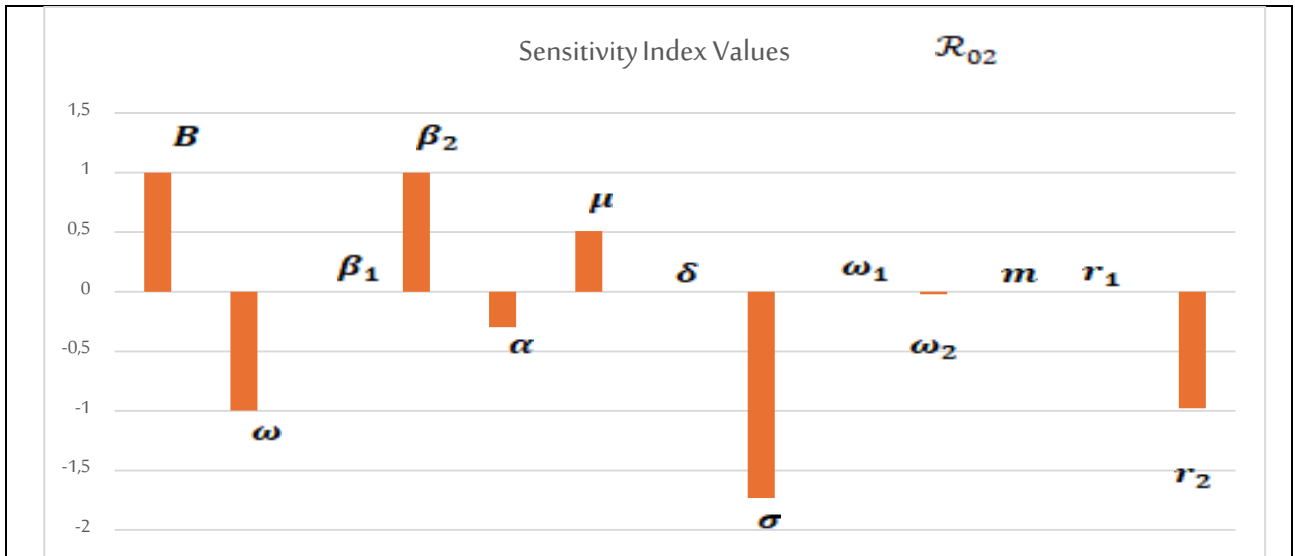

**Figure 7.2** Basic reproduction number $\mathcal{R}_{02}$ and sensitivity index values of the parameters

The parameters $B$, $\beta_1$, $\mu$ have positive sensitivity indices for $\mathcal{R}_{01}$, while the parameters $B, \beta_2, \mu$ have positive sensitivity indices for $\mathcal{R}_{02}$. Positive sensitivity indices indicate that an increase in the basic reproduction number is significant. In other words, increasing (or decreasing) the value of these parameters, while keeping other parameters constant, will lead to increases (or decreases) in the basic reproduction number. Conversely, the parameters $\omega, \alpha, \sigma, \omega_1, m, r_1$ have negative sensitivity indices for $\mathcal{R}_{01}$ and the parameters $\omega, \alpha, \sigma, \omega_2, r_2$ have negative sensitivity indices for $\mathcal{R}_{02}$. Negative sensitivity indices indicate that an increase in the basic reproduction number is of negative significance. In other words, increasing (or decreasing) the value of these parameters, while keeping other parameters constant, will result in decreases (or increases) in the basic reproduction number.

In conclusion, increasing the values of the parameters $B, \beta_1$ and $\mu$ for $\mathcal{R}_{01}$ and $B, \beta_2$ and $\mu$ for $\mathcal{R}_{02}$ indicates that the average number of secondary cases produced by a primary case in a fully susceptible population will increase. Similarly, decreasing the values of the parameters $\omega, \alpha, \sigma, \omega_1, m$ and $r_1$ for $\mathcal{R}_{01}$ and $\omega, \alpha, \sigma, \omega_2$ and $r_2$ for $\mathcal{R}_{02}$ shows that the average number of secondary cases produced by a primary case in a fully susceptible population will also increase.

The parameter to which the basic reproduction numbers $\mathcal{R}_{01}$ and $\mathcal{R}_{02}$ are most sensitive is $\sigma$ (vaccine efficacy). A reduction in the value of this parameter indicates that the average number of secondary cases transmitted by an infected individual in a fully susceptible population will increase. This suggests that as vaccine efficacy decreases, an infected person may transmit the disease to multiple individuals, leading to a progressive spread of the disease within the population. Our model highlights the importance of vaccine efficacy in epidemic control policies due to this negative sensitivity relationship with the basic reproduction number.

## 8 Simulations

In this section, we present various graphical representations that illustrate the dynamic behavior of the model system to demonstrate the underlying principles of the model dynamics and support the analytical results. We have considered the following initial conditions:

$$S(0) = 26195740, V = 51202223, I_1(0) = 269725, I_2(0) = 2724, R(0) = 7009861$$

**Simulation 1**

Here, we discuss the graphical representations of the behavior of the compartments infected with the original virus ($I_1$) and those infected with the mutant virus ($I_2$) for different values of the vaccine efficacy parameter ($\sigma$), along with an analysis of their impact on the epidemic.

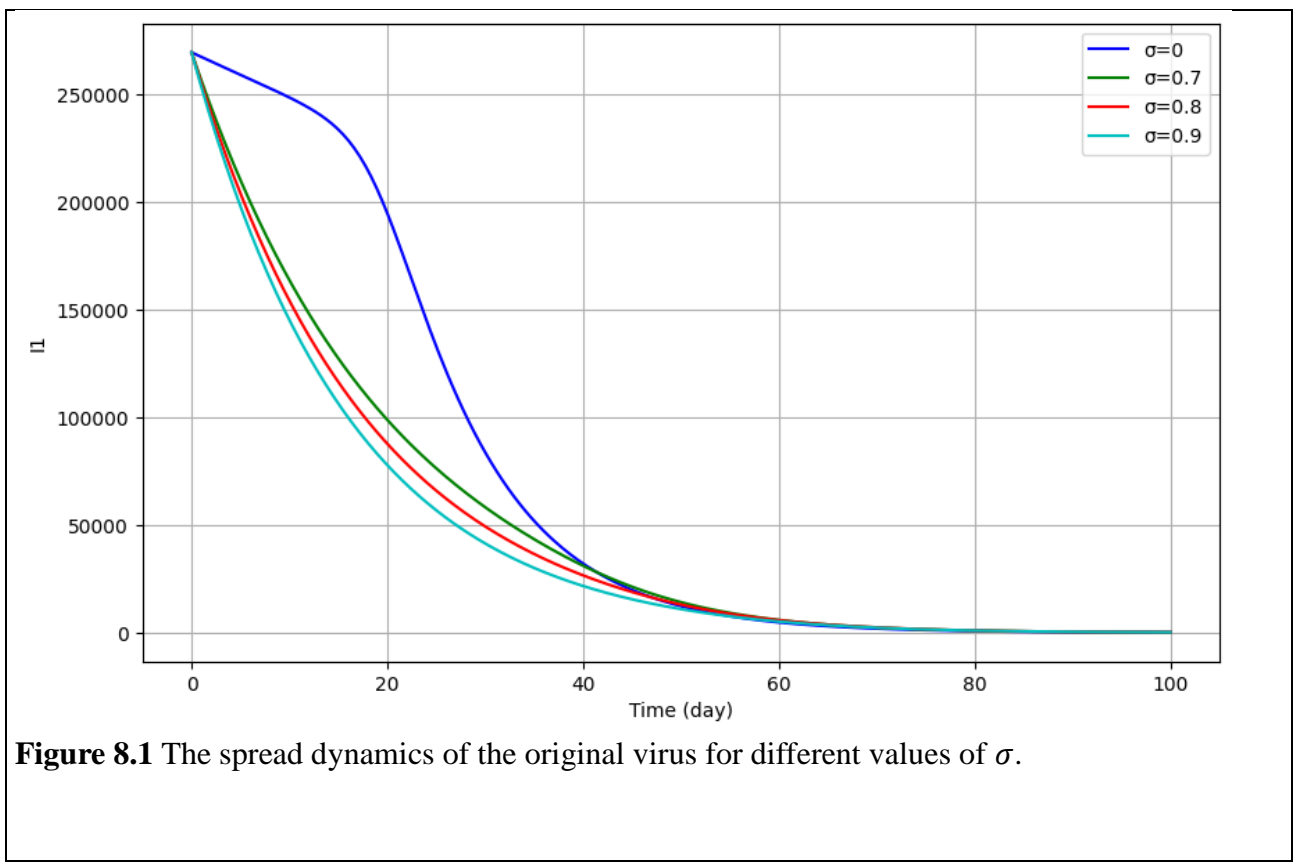

**Figure 8.1** The spread dynamics of the original virus for different values of $\sigma$.

In Figure 8.1, we examined the temporal variation of $I_1$ for different values of the vaccine efficacy parameter $\sigma$. The figure illustrates how the dynamics of individuals infected with the original virus evolve as the vaccine efficacy takes on different values. This analysis is crucial for understanding the role of vaccine efficacy in controlling the outbreak of the original virus and its impacts on the population.

All curves start from the same point on day 0, indicating a large number of infected individuals in the population who were either unvaccinated or unaffected by the vaccine at the outset. This reveals that the virus spread rapidly at the beginning of the outbreak, leading to a swift increase in the number of infected individuals.

$\sigma = 0$ represents the case where vaccine efficacy is zero, depicting a scenario in which the population is completely unvaccinated. The infection persists actively within the community for approximately 50 days before it starts to decline, although its complete eradication takes a longer period. This scenario demonstrates that in the absence of a vaccine, the virus continues to affect the population for an extended duration. In this case, the number of infected individuals reaches its peak and declines more slowly. This indicates that the outbreak will last longer and infect more individuals compared to scenarios where vaccines are present.

At $\sigma = 0.7$, although there is a reduction in the number of infected individuals, this decrease occurs more gradually, and the outbreak persists for a longer duration within the population. Around day 40, the number of infected individuals approaches zero.

At $\sigma = 0.8$, as vaccine efficacy increases, the reduction in the number of infected individuals accelerates further. With this value, the infection is observed to extinguish in a shorter period. Approximately 30 days later, the number of infected individuals approaches zero. This scenario underscores the significance of implementing a highly effective vaccine in controlling the spread of the virus within the community.

At $\sigma = 0.9$, there is a rapid decrease in the number of infected individuals, and the outbreak is controlled in a shorter time frame, reaching near zero in approximately 25 days. In this scenario, the vaccination program provides protection to a significant portion of the community, leading to a swift eradication of the virus. The rate of infection spread has nearly come to a standstill.

The graph serves as a clear demonstrator of the relationship between vaccine efficacy and the number of infected individuals, highlighting the critical importance of widespread vaccine use and efficacy for the future of epidemics. High vaccine efficacy (high values of σ) will lead to a rapid decrease in the number of infected individuals within the community, resulting in a shorter duration of the outbreak. This is significant for alleviating the burden on hospitals, minimizing economic impacts, and protecting overall public health. Such an effective vaccine also plays a crucial role in preventing new waves of outbreaks. Conversely, low vaccine efficacy (low values of σ) indicates that the number of infected individuals declines more slowly, leading to a longer duration for the outbreak to be controlled. This can result in prolonged infection and a greater number of affected individuals. Scenarios with low vaccine efficacy may also signal the emergence of new variants and the risk of resurgence of the outbreak.

In conclusion, there is a strong relationship between vaccine efficacy and the course of the outbreak. A more effective vaccination strategy plays a vital role in both shortening the duration of the outbreak and protecting public health. When vaccine efficacy is low, the likelihood of infected individuals remaining under the influence of the outbreak for an extended period increases, and the risk of long-term effects from the outbreak in the future grows.

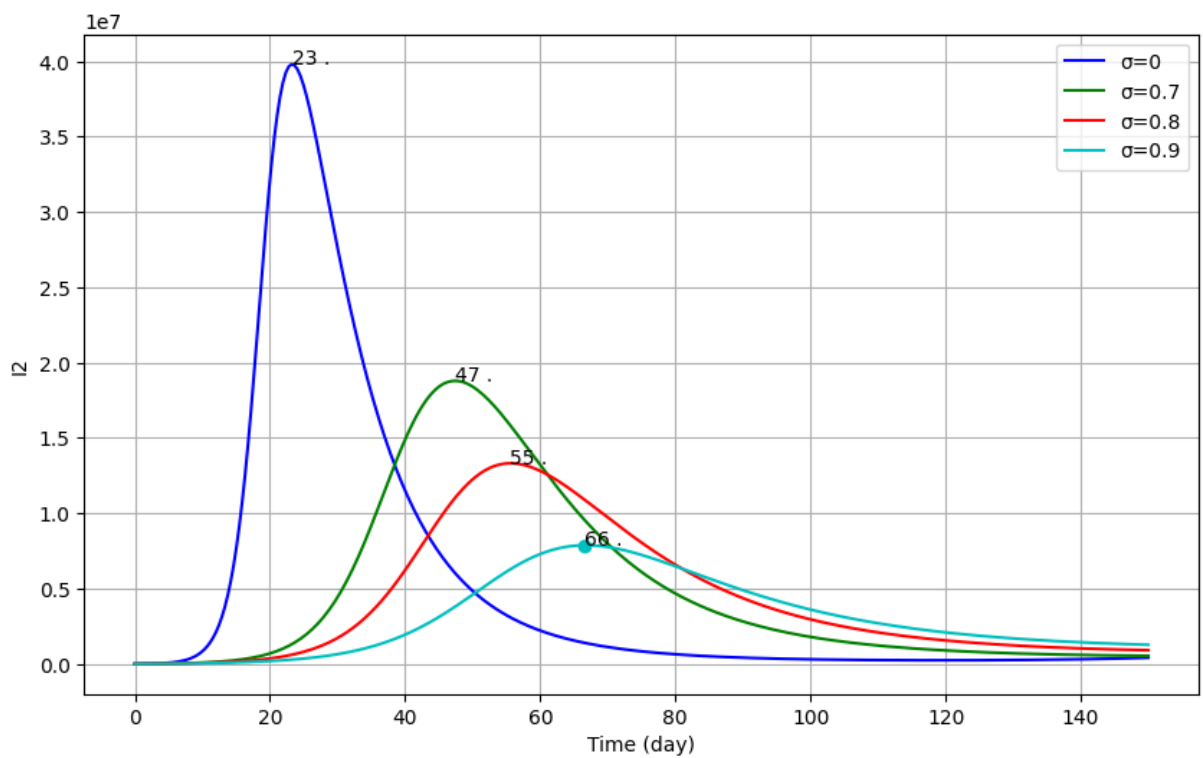

**Figure 8.2** The spread dynamics of the mutant virus for different values of σ.

Figure 8.2 shows the trend in the number of individuals infected with the mutant virus. In the absence of vaccination, $I_2$ reaches its peak on day 23, with approximately $4 \times 10^7$ individuals infected with the mutant virus. In the case of vaccination ($\sigma = 0.7$), the peak of $I_2$ occurs on day 47, with approximately $1.9 \times 10^7$ individuals infected. It can be observed that vaccination delays and reduces the peak of $I_2$. Thus, the peak of $I_2$ is postponed by 24 days. Even with a vaccine efficacy of $\sigma = 0.7$, the burden on medical resources caused by the outbreak will be significantly reduced, allowing more individuals to receive timely medical treatment, thereby decreasing mortality rates.

**Simulation 2**

Here, we discuss the graphical representations of the behavior of the compartments infected with the original virus $I_1$ and those infected with the mutant virus $I_2$ for different values of the vaccination rate parameter ($\alpha$), along with an analysis of their impact on the epidemic.

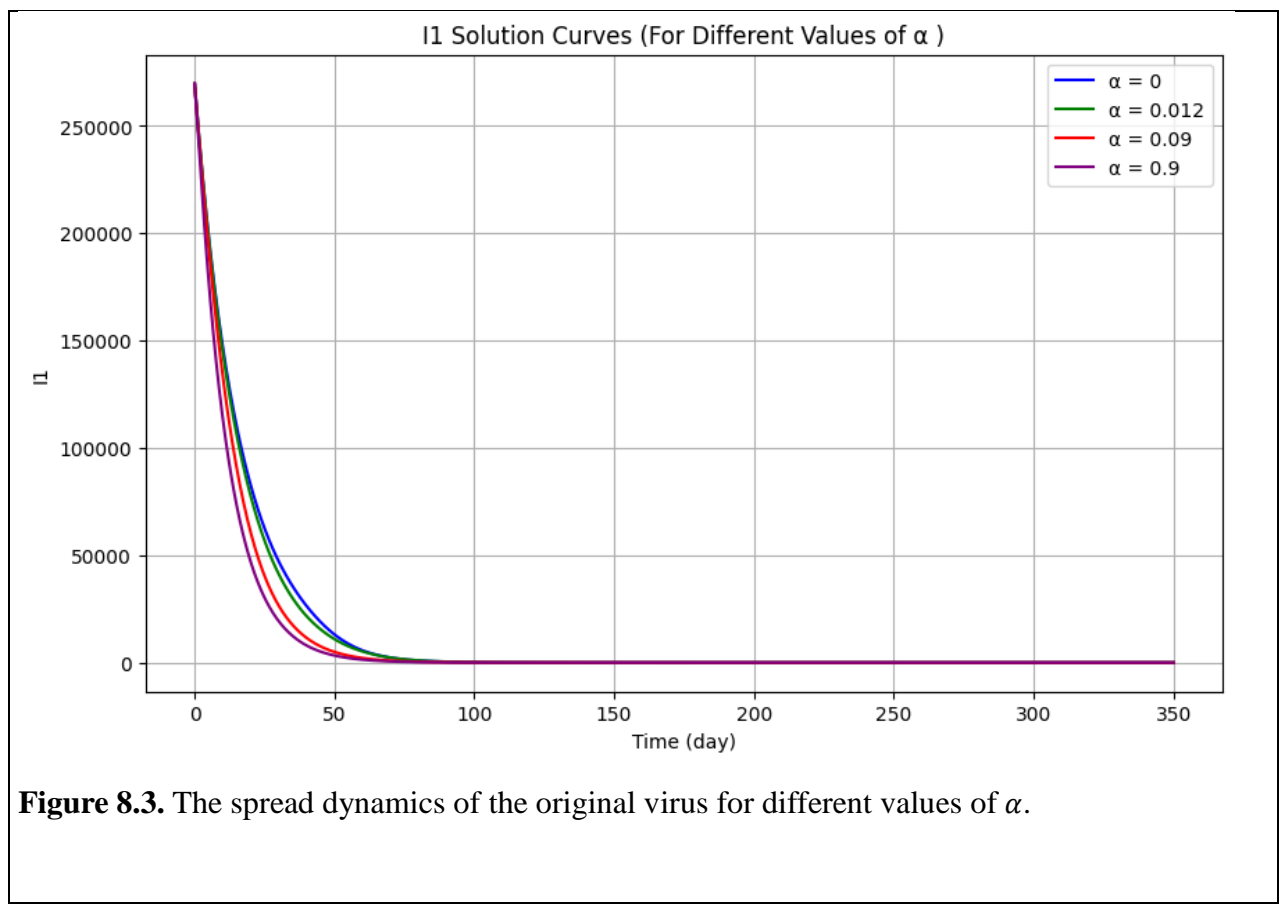

**Figure 8.3.** The spread dynamics of the original virus for different values of $\alpha$.

Figure 8.3 illustrates the temporal variation of $I_1$ for different values of the vaccination rate parameter $\alpha$. This graph clearly demonstrates the impact of vaccination on the reduction of individuals infected with the original virus over time. The higher the vaccination rate, the more rapidly individuals infected with the original virus are removed from the system. This highlights the significant role of vaccination in public health and its effectiveness in controlling the infection.

All curves in the graph start with a high number of infected individuals and decrease over time. In each case, the number of infected individuals approaches zero around day 100. This indicates that the infected individuals are either recovering or dying, thereby no longer remaining in the system.

At $\alpha = 0$, the number of individuals infected with the original virus shows a slower decline without vaccination and diverges somewhat from the other curves in the graph. This indicates that, in the scenario where vaccination is not implemented, individuals infected with the original virus remain in the system for a longer period.

At $\alpha = 0.012$, $\alpha = 0.09$ and $\alpha = 0.9$ the three curves remain closely aligned, indicating that as the vaccination rate increases, the number of individuals infected with the original virus declines more rapidly. High vaccination rates, such as $\alpha = 0.9$, lead to the fastest removal of individuals infected with the original virus from the system. However, even low and moderate vaccination rates contribute effectively to the reduction in the number of infected individuals.

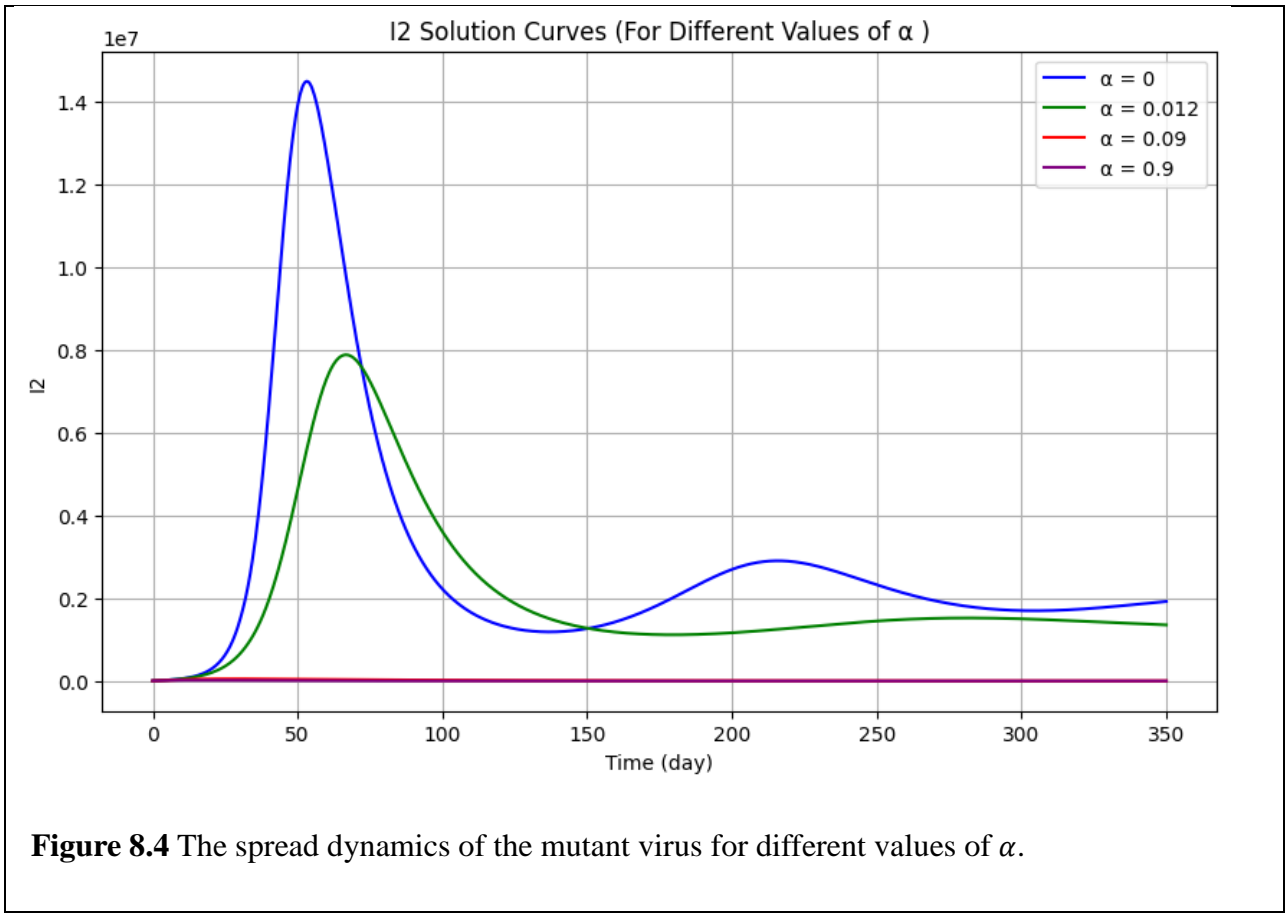

**Figure 8.4** The spread dynamics of the mutant virus for different values of $\alpha$.

The graph illustrates the temporal variation of the $I_2$ solution curves for different values of $\alpha$. The horizontal axis represents time (in days), while the vertical axis represents the values of $I_2$.

For $\alpha = 0$, there is an initial rapid increase that peaks around day 50, reaching a peak value of approximately $1.4 \times 10^7$. After this peak, there is a rapid decline; however, as time progresses, the decline slows down, and a steady state is observed around day 150. After day 150, a slight increase occurs again, and this fluctuation continues until the end of the graph.

For $\alpha = 0.012$, there is an initial rapid increase, but the peak value is lower (around $0.6 \times 10^7$) and occurs slightly earlier (around day 50). After the peak, the decline is less steep. Additionally, the fluctuations are less pronounced and have a lower amplitude.

For $\alpha = 0.09$, the curve remains significantly lower compared to the other two curves and continues almost horizontally. This indicates that as the value of $\alpha$ increases, the temporal variation of $I_2$ decreases substantially, leading to a more stable system.

For $\alpha = 0.9$, the curve remains almost at zero and shows no significant increase or decrease over time. This indicates that the system has reached complete equilibrium or remains within a very low range of variation.

In conclusion, the graph illustrates that as the value of $\alpha$ increases, there is a reduction in the temporal variation of $I_2$ and the system becomes more stable. Low $\alpha$ values (especially $\alpha = 0$) lead to larger amplitude and more fluctuating changes, whereas higher $\alpha$ values enhance the stability of the system and keep the variation at a minimum level.

**Simulation 3**

In this section, we explored the graphical representations and commentary on the behavior of the compartments infected with the original virus $I_1$ and the mutant virus $I_2$ for different values of the transmission rate of the original virus ($\beta_1$) and the transmission rate of the mutant virus ($\beta_2$), as well as their impact on the epidemic.

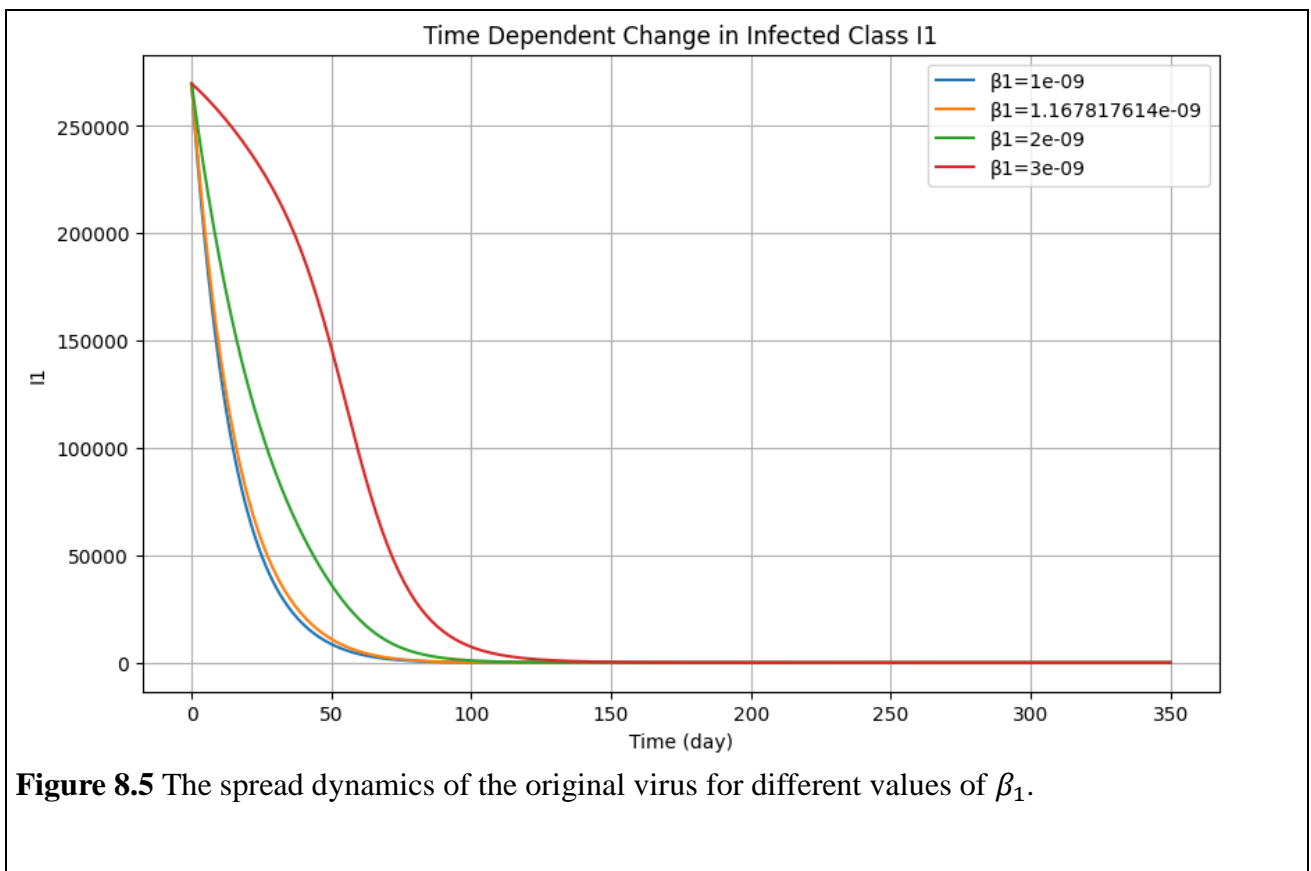

**Figure 8.5** The spread dynamics of the original virus for different values of $\beta_1$.

In Figure 8.5, the numerical variation of individuals infected with the original virus $I_1$ over time is examined under different values of the transmission rate $\beta_1$. The horizontal axis represents time (in days), while the vertical axis represents $I_1$, i.e., the number of infected individuals.

The graph shows that the number of individuals infected with the original virus will decrease over time, ultimately leading to the end of the outbreak. However, increases in the transmission rate significantly affect the duration and severity of the outbreak. At low transmission rates, the outbreak can be controlled more rapidly, while at high transmission rates, the effects of the outbreak persist longer and affect more individuals.

In the case of $\beta_1 = 0.000000001$, the transmission rate of the virus is extremely low. This results in a rapid initial decline in the number of infected individuals, reaching almost zero in approximately 75 days. This indicates that the outbreak can be quickly controlled at low transmission rates.

In the case of $\beta_1 = 0.000000001167817614$, the transmission rate is slightly higher, and the reduction in the number of infected individuals occurs more slowly. By approximately day 100, the number of infected individuals approaches zero. It is observed that as the transmission rate increases, the speed of virus spread also increases. Consequently, controlling the outbreak will take a bit longer in this scenario.

With a value of $\beta_1 = 0.000000002$, the transmission rate leads to a wider spread of the outbreak and causes infected individuals to remain in the system for a longer period. The number of infected individuals continues beyond day 100, and the rate of decrease occurs relatively slowly. In this scenario, the higher transmission rate of the virus indicates that the outbreak will persist longer and affect more individuals.

In the scenario with $\beta_1 = 0.000000003$, which has the highest transmission rate, the decline in the number of infected individuals occurs the slowest. The number of infected individuals remains relatively high initially and does not approach zero until approximately day 150. This indicates that in a scenario with a high transmission rate, the outbreak will persist longer, resulting in more individuals remaining infected and making control of the outbreak more challenging.

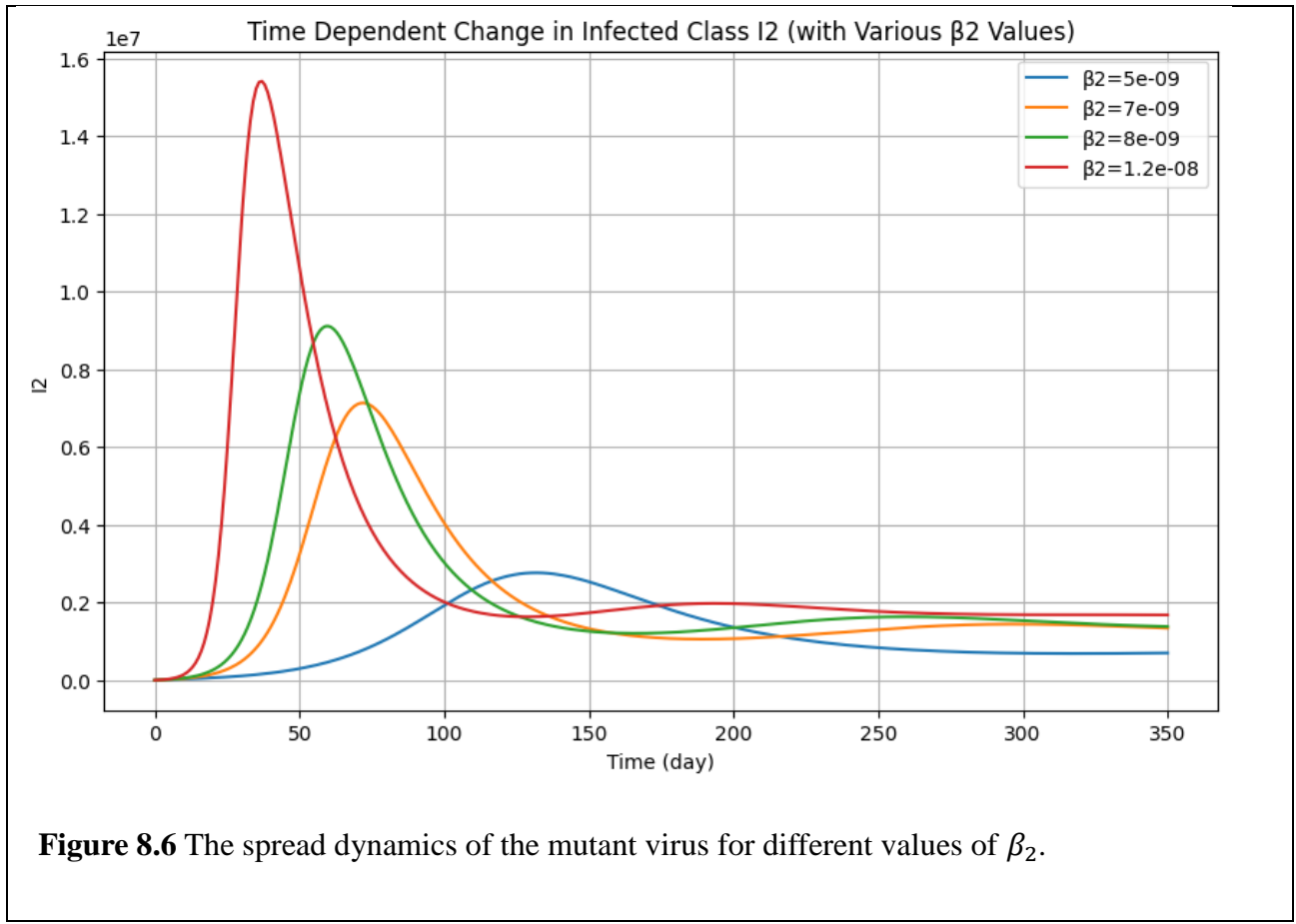

**Figure 8.6** The spread dynamics of the mutant virus for different values of $\beta_2$.

In Figure 8.6, we examined the numerical variation of individuals infected with the mutant virus ($I_2$) over time under different values of the transmission rate $\beta_2$. The horizontal axis represents time (in days), while the vertical axis expresses the number of infected individuals ($I_2$) on a logarithmic scale. This graph illustrates the spread dynamics of the mutant virus and how variations in the transmission rate affect the course of the outbreak.

The graph clearly demonstrates that the increase in the transmission rate of the mutant virus is one of the most critical factors determining the peak, duration, and severity of the outbreak. At lower transmission rates, the outbreak progresses more mildly and lasts longer. However, as the transmission rate increases, the outbreak spreads much more rapidly, the number of infected individuals reaches a higher peak, and there is significant potential to impact a large portion of the population in a short period, thereby placing a heavy burden on healthcare systems.

In the scenario with $\beta_2 = 0.000000005$, the lowest transmission rate, the number of infected individuals shows a very limited increase. It peaks at approximately day 75, with around $0.2 \times 10^7$ individuals infected. This transmission rate indicates that the spread of the mutant virus will be milder, resulting in fewer individuals becoming infected. Additionally, the duration of the outbreak is longer compared to other scenarios, as the lower transmission rate spreads the infected individuals over an extended period.

At a value of $\beta_2 = 0.000000007$, the number of infected individuals rises more rapidly, reaching its peak around day 60, with approximately $0.8 \, x \, 10^7$ (8 million) individuals infected. This increase in $\beta_2$ leads to a more severe wave of the outbreak. The higher number of infected individuals and their shorter duration in the system indicate that the virus will spread quickly, potentially affecting more individuals in a short period.

At a value of $\beta_2 = 0.000000008$, the peak number of infected individuals occurs around day 50. In this scenario, the peak value is approximately $1 \, x \, 10^7$ (10 million) individuals, further increasing the severity of the outbreak. With this increase, it is observed that infected individuals multiply more rapidly and the outbreak reaches its peak in a shorter time. The speed of the outbreak's spread, driven by the increase in the transmission rate, results in more individuals becoming infected and places greater strain on healthcare systems.

At a value of $\beta_2 = 0.000000012$, the highest transmission rate, the number of infected individuals increases most rapidly, reaching its peak around day 45. At this transmission rate, it is evident that the mutant virus spreads quickly and affects a large number of individuals. However, after peaking, the decrease in the number of infected individuals is also quite rapid. This indicates that while the virus is highly contagious, a significant portion of the population becomes infected in a short time, leading to a quicker decline of the outbreak.

## 9 Conclusions

In this study, the deterministic compartmental model we developed to understand the dynamics of COVID-19 infection provides a comprehensive framework for epidemiological analyses. Our model consists of compartments for $S$ (susceptible), $V$ (vaccinated), $I_1$ (infected with the original virus), $I_2$ (infected with the mutant virus), and $R$ (recovered), allowing us to examine the spread of the virus in detail.

The basic reproduction numbers obtained, $R_{01}$ and $R_{02}$, play a critical role in determining the potential spread of the original and mutant viruses within the population. When $R_{01}$ and $R_{02}$ are less than one, the global asymptotic stability of the disease-free equilibrium indicates that the viruses can be controlled, whereas values greater than one signify the stability of the endemic equilibrium. This suggests that the viruses can sustain their transmission within the community, highlighting the urgent need for public health policy reassessment.

Bifurcation analysis demonstrates the presence of backward bifurcation for the original virus, shedding light on the complexity of epidemiological dynamics and the impact of parameter changes. Notably, increases in transmissibility rates lead to rapid rises in the number of infected individuals in the $I_1$ and $I_2$ compartments. This increases the risk of more individuals becoming infected as the viruses spread through the community, significantly elevating the burden on health systems. High transmissibility rates not only affect the number of infected individuals but also the speed and duration of disease spread, posing a critical threat to public health.

Furthermore, the results from the sensitivity analysis of the vaccine efficacy ($\sigma$) parameter emphasize the importance of vaccination in public health policies. When vaccine efficacy is high, a notable decrease in the number of infected individuals is observed, enabling better control of the outbreak. Conversely, low vaccine efficacy results in infected individuals remaining in the system for longer periods, leading to wider community impact. This highlights the necessity for optimizing vaccination strategies and increasing vaccination rates across the population.

Our study provides a detailed analysis of various parameters (such as vaccine efficacy and transmissibility rates) that affect the course of the COVID-19 pandemic in Turkey. The validation of the model with real data from COVID-19 cases in Turkey enhances the validity of the findings and establishes a solid foundation for future epidemiological studies.

In conclusion, this work contributes to a better understanding of the dynamics of COVID-19 through the developed model and aids in optimizing public health strategies. Future research should comprehensively investigate the effects of mutant viruses and vaccines, contributing to the development of more effective public health measures. Additionally, testing and updating the model under different scenarios will be a significant step in preparing for future pandemics. Thus, the findings provide critical information for decision-makers in the fight against COVID-19 and illuminate efforts to protect public health.

**Appendix A**

$k_1 = -d_{11} - d_{22} - d_{33} - d_{44} - d_{55},$

$k_2 = d_{11}d_{22} - d_{12}d_{21} - d_{31}d_{13} - d_{14}d_{41} - d_{23}d_{32} - d_{24}d_{42} + d_{44}(d_{11} + d_{22} + d_{33}) - d_{33}(-d_{11} - d_{22}) - d_{55}(-d_{11} - d_{22} - d_{33} - d_{44}),$

$k_3 = -d_{15}d_{53}d_{31} - d_{15}d_{54}d_{41} - d_{31}(d_{11}d_{13} + d_{12}d_{23}) - d_{32}(d_{21}d_{13} + d_{22}d_{23}) - d_{41}(d_{11}d_{14} + d_{12}d_{24}) - d_{42}(d_{21}d_{14} + d_{22}d_{24}) - d_{43}(d_{31}d_{14} + d_{32}d_{24}) + (d_{11} + d_{22} + d_{33})(d_{14}d_{41} + d_{24}d_{42}) - d_{33}(d_{11}d_{22} - d_{12}d_{21}) - (d_{13}d_{31} + d_{23}d_{32})(-d_{11} - d_{22}) + d_{44}(-d_{11}d_{22} + d_{12}d_{21} + d_{13}d_{31} + d_{23}d_{32} + d_{33}(-d_{11} - d_{22})) - d_{55}(d_{11}d_{22} - d_{12}d_{21} - d_{13}d_{31} - d_{14}d_{41} - d_{23}d_{32} - d_{24}d_{42} + d_{44}(+d_{11} + d_{22} + d_{33}) - d_{33}(-d_{11} - d_{22})),$

$k_4 = -d_{53}(d_{11}d_{31}d_{15} + d_{21}d_{32}d_{15} + d_{31}d_{15}d_{33}) - d_{54}(d_{11}d_{41}d_{15} + d_{21}d_{15}d_{42} + d_{31}d_{15}d_{43} - d_{41}d_{15}d_{44}) - d_{41}(d_{11}(d_{11}d_{14} + d_{12}d_{24}) + d_{12}(d_{21}d_{14} + d_{22}d_{24}) + d_{13}(d_{31}d_{14} + d_{32}d_{24})) - d_{42}(d_{21}(d_{11}d_{14} + d_{12}d_{24}) + d_{22}(d_{21}d_{14} + d_{22}d_{24}) + d_{23}(d_{31}d_{14} + d_{32}d_{24})) - d_{43}(d_{31}(d_{11}d_{14} + d_{12}d_{24}) + d_{32}(d_{21}d_{14} + d_{22}d_{24}) + d_{33}(d_{31}d_{14} + d_{32}d_{24})) + (d_{11} + d_{22} + d_{33})(d_{41}(d_{11}d_{14} + d_{12}d_{24}) + d_{42}(d_{21}d_{14} + d_{22}d_{24}) + d_{43}(d_{31}d_{14} + d_{32}d_{24})) + (d_{14}d_{41} + d_{24}d_{42})(-d_{11}d_{22} + d_{12}d_{21} + d_{13}d_{31} + d_{23}d_{32} + d_{33}(-d_{11} - d_{22})) + d_{44}(d_{31}(d_{11}d_{13} + d_{12}d_{23}) + d_{32}(d_{21}d_{13} + d_{22}d_{23}) + d_{33}(d_{11}d_{22} - d_{12}d_{21}) + (d_{13}d_{31} + d_{23}d_{32})(-d_{11} - d_{22})) - (d_{31}d_{15}d_{53} + d_{41}d_{15}d_{54})(-d_{11} - d_{22} - d_{33} - d_{44}) - d_{55}\big(-d_{31}(d_{11}d_{13} + d_{12}d_{23}) - d_{32}(d_{21}d_{13} + d_{22}d_{23}) - d_{41}(d_{11}d_{14} + d_{12}d_{24}) - d_{42}(d_{21}d_{14} + d_{22}d_{24}) - d_{43}(d_{31}d_{14} + d_{32}d_{24}) + (d_{11} + d_{22} + d_{33})(d_{14}d_{41} + d_{24}d_{42}) - d_{33}(d_{11}d_{22} - d_{12}d_{21}) - (d_{13}d_{31} + d_{23}d_{32})(-d_{11} - d_{22}) + d_{44}(-d_{11}d_{22} + d_{12}d_{21} + d_{13}d_{31} + d_{23}d_{32} + d_{33}(-d_{11} - d_{22}))\big),$

**Acknowledgements**

Not applicable

**Funding**

Not applicable

**Availability of data and materials**

All data are publicly available in the referenced sources.

**Declarations**

**Competing interests**

The authors declare that they have no competing interests.

**Authors' contributions**

CG conceived the study, CG and SŞ carried out the analysis, CG wrote the draft, all authors revised the manuscript critically, and approved it for publishing.